\documentclass[11pt]{article}

\usepackage{amsmath,amsthm}
\usepackage{amssymb}
\usepackage[mathscr]{eucal}
\usepackage[all]{xy}
\usepackage{mathrsfs}
\usepackage{hyperref}
\usepackage{appendix}
\usepackage{authblk}
\usepackage{graphicx}

\usepackage{hyperref}

\usepackage[text={160mm,240mm},centering]{geometry}

\newtheorem{theorem}{Theorem}

\newtheorem{proposition}[theorem]{Proposition}

\theoremstyle{definition}

\newtheorem{definition}[theorem]{Definition}
\newtheorem{definition-lemma}[theorem]{Definition-Lemma}
\newtheorem{definition-theorem}[theorem]{Definition-Theorem}

\newtheorem{remark}[theorem]{Remark}
\newtheorem*{remarks}{Remark}

\newtheorem*{ack}{Acknowledgements}

\def\be{\begin{equation}}
\def\ee{\end{equation}}
\def\bea{\begin{eqnarray}}
\def\eea{\end{eqnarray}}
\def\bes{\begin{eqnarray*}}
\def\ees{\end{eqnarray*}}

\begin{document}

\title{Gutzwiller's Semiclassical Trace Formula \\ and \\ Maslov-Type Index Theory for Symplectic Paths}

\author[1,2]{Shanzhong Sun\thanks{Partially supported by NSFC (No.10731080, 11131004, 11271269),
PHR201106118, the Institute of Mathematics and Interdisciplinary Science at CNU, Email: sunsz@cnu.edu.cn}}

\renewcommand\Affilfont{\small}

\affil[1]{Department of Mathematics, Capital Normal University, Beijing 100048 P. R. China}
\affil[2]{Beijing Center for Mathematics and Information Interdisciplinary Sciences, Beijing 100048 P. R. China}

\date{}

\maketitle
\begin{center}
{Dedicated to Professor Paul Henry Rabinowitz with admiration}
\end{center}

\begin{abstract}
Gutzwiller's famous semiclassical trace formula plays an important role in theoretical and experimental quantum mechanics with tremendous success. We review the physical derivation of this deep periodic orbit theory in terms of the phase space formulation with an view towards the Hamiltonian dynamical systems. The Maslov phase appearing in the trace formula is clarified by Meinrenken as Conley-Zhender index for periodic orbits of Hamiltonian systems. We also survey and compare various versions of Maslov indices to establish this fact. A refinement and improvement to Conley-Zehnder's index theory which we will recall all essential ingredients is the Maslov-type index theory for symplectic paths developed by Long and his collaborators which would shed new light on the computations and understandings on the semiclassical trace formula. The insights in Gutzwiller's work also seems plausible to the studies on Hamiltonian systems.
\end{abstract}

\setcounter{tocdepth}{2}\tableofcontents

\section{Introduction}
\renewcommand{\thetheorem}{\Alph{theorem}}

Classical mechanics has evolved from the Newtonian formulation to Lagrangian variational formulation and Hamiltonian symplectic geometric formulation since the inception of the science from the late 17th century. Special roles are played by periodic orbits which are our main concern. In his studies on three-body problem in celestial mechanics, H. Poincar\'{e} not only discovered the phenomenon of (classical) chaos, a term encoding the extremely irregular character of the orbit, but also raised the less famous (compared with his conjecture on three-sphere in topology which is proved recently by G. Perelman) Poincar\'{e} conjecture on the denseness of periodic orbits in the restricted three-body problem. The general importance of periodic orbits was emphasized in his three-volume monumental monograph "The New Methods of Celestial Mechanics". Other than very few systems showing regular and simple behavior, we have to resort to perturbation theory. The famous example is the powerful Kolmogorov-Arnold-Moser theory which claims the survival of regular behavior under some small perturbation, and can be used to detect the transition from the regular system to chaotic behavior. Now we know that even a simple dynamical system can show chaotic irregularity, and chaotic systems lie out of reach of perturbative analysis. Examples of chaotic classical systems are double pendulum, billiard, geodesics on hyperbolic surface and planetary dynamics among others.

At about the same time of Poincar\'{e}, M. Plank, A. Einstein and N. Bohr revolutionized physics by the idea of quantum. Modern theory of quantum mechanics was firmly established in the pioneering works of L. de Broglie, W. Heisenberg, E. Schr\"{o}dinger, M. Born, W. Pauli and P. Dirac in five years from 1924 to 1929. It is a field full of incredible subtleties and difficulties. Most of twentieth century mathematics are developed for this purpose, say functional analysis and operator algebra, representation theory and noncommutative geometry to name a few.

It is always an intriguing problem to understand the relationship between classical mechanics and quantum mechanics. It is also absolutely needed to understand modern physical adventures like supersymmetric gauge theory and mirror symmetry in string theory.

Combining well known classical mechanical laws with correspondence principle, Bohr succeeded in establishing the quantization of hydrogen atom and getting the quantum mechanical energy spectra of the Coulomb potential. It is not only of historical importance in motivating the adventures of quantum mechanics, but also inspired WKB quantization conditions for classical integrable systems of several degree of freedoms with culmination in Maslov's work on phase loss when bypassing the caustics.

Already in 1917 at the beginning of the development of the quantum mechanics, Einstein (\cite{E17}) had raised the question on the relation between classically chaotic systems and quantum mechanics. More precisely, is there any quantum analogue of classical chaos?  He clearly observed that the Bohr correspondence principle could not be directly extended to the chaotic system. One of the prominent feature of quantum mechanics is the wave like nature of the atom scale systems. It deems to be a difficult problem to reconcile these two seemingly quite different phenomenon. It was forgotten for more than forty years. In 1970's, M. Gutzwiller was the first to study the relationship between classical and quantum mechanics in chaotic systems. All the attempt to understand this transition back and forth is called quantum chaos/quantum chaology which tries to build a bridge between classical mechanics and quantum mechanics. Main topics are distribution of quantized energy levels, stationary states. A current tantalizing hot topic is the mysterious and deep parallel between the chaotic scattering in quantum chaos and the Riemann hypothesis: the distribution of zeros of the Riemann's zeta function in number theory (see e.g. \cite{BK13} for the state of the art and references therein). Physical applications include the energy levels of the donor impurity in a silicon (Si) or germanium (Ge) crystal (anisotropic Kepler problem, \cite{G73}), ordinary hydrogen atom near ionization in a strong magnetic field with the discovery of the Resurgence Spectrocopy; Rydberg hydrogen atom;  statistical properties of the spectra of nucleus. Note that random matrix theory was introduced into this field by M. Berry generalizing the idea of E. P. Wigner and F. J. Dyson which is another main tools of theoretical investigations of quantum chaotic systems. In this latter aspect, they were among the first people to study the universal statistical properties of nuclear spectra of a chaotic quantum system (i.e. quantum system with chaotic classical analogue), and they found the most probable formula for the distribution as confirmed by Bohigas and Giannoni which we will not touch upon in this work.

To get the semiclassical spectra of a chaotic quantum system,  Gutzwiller(\cite{G71}) found a very general, beautiful and deep trace formula which extracts the information on the eigenvalues of the Hamiltonian operator in terms of the complete enumeration of the periodic orbits of its corresponding classical Hamiltonian system. The moral is that the classical periodic orbits and the quantum mechanical spectra are closely tied together through Fourier transform. This put the periodic solutions of classical Hamiltonian system as emphasized by Poincar\'{e} into a new perspective. It is called in the literature Gutzwiller's semicalssical trace formula (SCTF) or his periodic orbit theory (POT). It has been successfully used by him and Berry among others to study the statistical properties of the spectrum and low-lying eigenvalues. We should point that the SCTF plays a mutual role, in other words, it can also be used in the reverse order to get the information on periodic orbits via energy spectrum (inverse quantum chaology/quantum  recurrence spectroscopy). An early example is due to J. Chazarain (\cite{C74}) who showed that summation over energy levels of the Laplace operator generated a function with singularities the actions of closed geodesics. For more recent progress in this direction, please refer to \cite{KS14} and references therein.

For a Hamiltonian system with Hamiltonian $H$, suppose that the corresponding operator $\hat H_{\hbar}$ has a discrete eigenvalues sequence $E_0,E_1,E_2,...,E_n,...$. Roughly speaking, Gutzwiller's formula says that when $\hbar\rightarrow 0$, the density of states $\rho(E):=\sum_{j=0}^\infty\delta(E-E_j)$ is approximated, up to the so-called "Weyl term" $\langle\rho(E)\rangle$ which counts the number of states in the phase space region $H(z)\leq E$, by the sum $$\frac{1}{\pi\hbar}\textrm{Re}\sum_{\gamma}\sqrt{-1}^{-i_\gamma}\frac{T_\gamma e^{\frac{\sqrt{-1}}{\hbar}\mathcal{A}_\gamma}}{{\sqrt{|\textrm{det}(P_\gamma-I)|}}},$$ where $\gamma$ runs over all closed periodic orbits of $H$ with energy $E$ and period $T_\gamma$ (including their iterations). We assume that this set consists only of non-degenerate orbits, i.e. the linearized Poincar\'{e} map $P_\gamma$ for $\gamma$ has no eigenvalue $1$.  $\mathcal{A}_\gamma=\oint_\gamma pdq$ is the classical action. $i_\gamma$ is the Conley-Zehnder index of a certain symplectic path $\tilde P_\gamma$ associated to $P_\gamma$(\cite{M92,M94}), which is our main concern in the current paper.

Heuristically, one of the physical derivations is through the alternative evaluation of the path integral in quantum mechanics. The path integral is an infinite dimensional integral over the path space with all possible end points in configuration space. One way is to replace time $t$ by energy $E$ via a Fourier transform, and taking trace is equivalent to take integrals over free loop spaces resulting in the energy spectrum density. On the other hand, one can perform the semiclassical approximation which receives contributions only from the closed periodic orbits in the corresponding classical system. Here we repeat Gutzwiller (\cite{G07}) to give an intuitive interpretation: "The open parameter $E$ represents a small perturbation with a constant frequency $\mu=E/h$ that works on the system from the outside, where $h$ is always Planck's constant. The reaction of the system is a forced motion of the same frequency, with the amplitude $g(E)$. The closure $E$ is to one of the eigenvalues $E_n$, the larger is the response of the system; we get a resonance! The external perturbation of frequency $\mu$ can be described also by its period $\tau$, the reciprocal of $\mu$. The classical particle gets chased around in its space, and it is critical where it lands after one period $\tau$. The effect on the classical particle will be larger if it comes back to its starting point after one, or perhaps two or three such periods. Therefore, the classical description of a quantum resonance depends on the periodic orbits. The physical length of a periodic orbit yields the period in time by taking the derivative with respect to the energy $E$ of the periodic orbit." There are several proofs to this SCTF with various mathematical rigorousness which we will comment in later sections.

Several special cases of the SCTF are already very interesting. Notably among them is the famous Selberg trace formula which is even an exact formula and not just an asymptotic expansion, a miracle! We will come back to this formula later on. In fact, A. Selberg got his formula for locally symmetric spaces. It was H. Huber (\cite{H59}) who interpreted it as a formula relating the spectrum of the Laplace operator and the lengths spectrum on a closed Riemann surface of constant curvature $-1$. This is the first instance of the SCTF in the history.

In physics, at about the same time as Gutzwiller, R. Balian and C. Bloch (\cite{BB70}, \cite{BB71}, \cite{BB72} and \cite{BB74}) developed also the asymptotic trace formula for the eigenfrequencies of a cavity ("multiple reflection expansion" in their terminology) independently. This is closely related to the famous M. Kac's problem: Can one hear the shape of a drum? namely the inverse spectral problem.

On the mathematical side, just after the publications of the works of the physicists (around 1973-1975), several groups of mathematicians devoted to the rigorous studies for Laplace operator on closed Riemannian manifolds. Colin de Verdi\`{e}re in his thesis (see \cite{CdV73a}, \cite{CdV73b}) proved that the spectrum of the Laplacian determines generically the lengths spectrum by using the short-time expansion of the Schr\"{o}dinger kernel and a discrete approximation of Feynman path integral. Chazarain (\cite{C74}) got the qualitative form of the trace for the wave kernel by using the Fourier integral operators. J. J. Duistermaat and V. Guillemin (\cite{DG75}) studied the main term of the singularity in terms of the Poincar\'{e} map of the closed orbit by using the full power of the Fourier integral operators, and now \cite{DG75} is a standard reference on the subject. Later on, it is extended to more general semiclassical Hamiltonians by Helffer-Robert, Guillemin-Uribe and Meirenken among others, to manifolds with boundaries, to several commuting operators by Charbonnel-Popov. Recently, wave trace theory is developed by Guillemin and S. Zelditch (\cite{G96}, \cite{Z97}, \cite{Z98}) which can be used to draw information on the nonprincipal terms in the singularities expansion by the semiclassical Birkhoff normal forms.

Some surveys related to the trace formula are \cite{CdV98}, \cite{CdV06}, \cite{CdV07}, \cite{G77}, \cite{G90}, \cite{G07}, \cite{U00}, and \cite{W}.

The appearance of Morse index and Maslov index has been noticed at the very beginning in the development of semiclassical quantization.  There are many works in physics devoted to the understanding and computations about Maslov index appearing in the semiclassical trace formula (e.g. \cite{CRL90}, \cite{LR87}, \cite{BP03}, \cite{R92}, \cite{RL87}, \cite{S99}, \cite{S02} to name a few).

In mathematics, there are two essentially equivalent ways to define the Maslov index. One is in the framework of Lagrangian Grassmanian which is introduced by V. P. Maslov(\cite{M72}) in 1965 in the multidimensional asymptotic analysis of semiclassical quantization. The other is in the framework of symplectic group initiated by I. M. Gelfand and V. B. Lidskii (\cite{GL55}, index of rotation) in 1955 in their studies of stability of linear periodic Hamiltonian systems; in 1984 C. Conley and E. Zehnder (\cite{CZ84}) put it in the variational perspective. V. I. Arnold (\cite{A67}) discovered the equivalence between these two different formulations in 1967. Some of the facts to establish this equivalence will be recalled later on in this paper. In the mathematical literature, there are several versions of Maslov index, for example, that of 1) a smooth one-parameter family of pairs of Lagrangian subspaces, 2) triple of Lagrangian  subspaces,  3) a path of symplectic matrices. Although they are all related to each other, the point here is to decide the right version for the problem at hand.

In this note, we do not claim any originality. The purpose of the paper is two-fold. On the one hand we try to put the research on periodic orbits of general Hamiltonian system into a new perspective or better an old tradition; and on the other hand we hope that the recent advances in classical Hamiltonian dynamics, especially the Maslov-type index theory developed by Long and his collaborators will be helpful for computations and understanding on the spectrum of Hamiltonian operators in quantum mechanics. In this sense, the paper can be seen as a guide to the excellent book of Long (\cite{Long02}). There is a huge literature devoted to topics as the Maslov indices and quantum chaos, and a complete list of references is almost impossible. We are content with the most relevant papers and monographs which are easy to access with further information on references.

The paper is organized as follows: After introduction to the general idea on Gutzwiller's semiclassical trace formula and Maslov index in \S 1, we review various versions of the definitions of Maslov index from Lagrangian perspective in \S 2 and symplectic group perspective in \S3 which include a detailed description on the Maslov-type index theory developed by Long and their relationships. In \S 4, the relation between Morse index and Maslov index is  given to get an idea about this deep phenomenon. In \S 5, we show the physical heuristic derivation of the semiclassical trace formula in phase space formulation which still calls for mathematically rigorous justification, several attempts among mathematicians via microlocal analysis is briefly summarized, a model example, i.e., the well-known Selberg trace formula is analyzed, and the role of Maslov-type index theory is clarified. We conclude the paper with several promising directions to study further this fascinating trace formula in \S 6.

\begin{ack}
I would like to thank Professor Yiming Long for teaching me the fine points in Maslov-type index theory for symplectic paths, and for many encouragements and supports during the past two decades. Thanks also go to the referee for his/her careful readings and valuable suggestions.

\end{ack}

\setcounter{theorem}{0}
\renewcommand{\thetheorem}{\arabic{theorem}}

\section{Maslov Index: Lagrangian Grassmannian Perspective}

For a real Symplectic vector space $(V,\omega)$, we denote by $\mathrm{Lag}(V,\omega)$ its Lagrangian Grassmannian and $\mathrm{Sp}(V,\omega)$ its symplectic group.

Since all symplectic vector spaces of the same dimension are isomorphic to each other by classical Darboux theorem, we use in our most descriptions the standard symplectic vector space $E:=(\mathbf{R}^{2n},\omega)$, denote its Lagrangian Grassmannian by $\Lambda(n)$ and the symplectic group by $$\mathrm{Sp}(2n):=\{M\in GL(2n,\mathbf{R})\,|\, M^\top JM=J\}$$ with $J=\left(\begin{matrix}0 & -I_n \cr I_n & 0\cr\end{matrix}\right)$ and $I_n$ the identity matrix on $\mathbf{R}^n$. The spectrum/eigenvalues of a symplectic matrix $M$ is $\sigma(M):=\{\lambda\in\mathbf{C}\,|\,\det(M-\lambda I_{2n})=0\}$.

\subsection{Lagrangian Grassmannian}

\cite{A67}, \cite{RS93}

We recall some topological facts about Lagrangian Grassmannian $\Lambda(n)$.

\begin{proposition}
The unitary group $U(n)$ acts on $\Lambda(n)$ transitively with stabilizer group $O(n)$.
\end{proposition}
So $\Lambda(n)$ is a manifold with homogeneous space structure $\Lambda(n)=U(n)/O(n)$ with the canonical map $\textrm{Det}^2: \Lambda(n)\rightarrow S^1$. Using homotopy exact sequence one can prove that

\begin{proposition} The first homology group, cohomology group and the fundamental group of $\Lambda(n)$ are the same:
 $$H_1(\Lambda(n),\mathbf{Z})\cong H^1(\Lambda(n),\mathbf{Z})\cong \pi_1(\Lambda(n))\cong\mathbf{Z}.$$
\end{proposition}

One can find the Maslov singular cycle and cocycle as the generators of the first (co)homology groups. {\bf Maslov index} of a closed curve $\gamma: S^1\rightarrow \Lambda(n)$ can be defined as its topological intersection number with singular cycle, or as the rotation number of $\textrm{Det}^2$, or as the mapping degree of $\textrm{Det}^2\circ \gamma: S^1\rightarrow S^1$.

Denote by $\Lambda^k(n)$ the subset of $\Lambda(n)$ consisting of Lagrange subspaces having a $k(\geq 0)$-dimensional intersection with a fixed $L_0\in\Lambda(n)$. $\Lambda^0(n)$ is an open set in $\Lambda(n)$ which is diffeomorphic to the linear space of all real symmetric matrices of order $n$. The set $\Lambda^k(n)$ is an open manifold of codimension $k(k+1)/2$ in $\Lambda(n)$. $\overline{\Lambda^1(n)}=\cup_{k\geq 1}\Lambda^k(n)$ gives the Maslov singular cycle which is two-sidely embedded in $\Lambda(n)$.

\subsection{Triple Maslov Index}

H\"{o}rmander(\cite{H71}), Kashiwara(\cite{KS90}(Appendix), \cite{LV80}), \cite{GS77}

For ordered Lagrangian triplets $(L_1,L_2,L_3)\in\Lambda(n)^3$, we define the quadratic form \begin{eqnarray*} Q(L_1,L_2,L_3): L_1\oplus L_2\oplus L_3&\rightarrow& \mathbf{R}\\
(x_1,x_2,x_3) &\mapsto& \omega(x_1,x_2)+\omega(x_2,x_3)+\omega(x_3,x_1).\end{eqnarray*}
The {\bf triple Maslov index (or inertia index)} is defined to be the signature of the quadratic form of the triplets of Lagrangian subspaces, namely,
$$ s(L_1,L_2,L_3):= \mathrm{sgn}(Q(L_1,L_2,L_3)).$$

\begin{proposition}\label{signature}
The signature $s$ has the following properties:
\begin{itemize}
\item[(1)] invariant under the action of $\mathrm{Sp}(2n)$: $$s(gL_1,gL_2,gL_3)=s(L_1,L_2,L_3), \forall g\in \mathrm{Sp}(2n),$$
  which means the diagonal action of $\mathrm{Sp}(2n)$ on $\Lambda(n)^3$ is not transitive, contrary to the transitivity of the action on $\Lambda(n)^2$;
\item[(2)] antisymmetric under the permutations of the Lagrangians: $$s(L_1,L_2,L_3)=-s(L_2,L_1,L_3)=-s(L_1,L_3,L_2);$$
\item[(3)] Cocycle identity: $s(L_2,L_3,L_4)-s(L_1,L_3,L_4)+s(L_1,L_2,L_4)-s(L_1,L_2,L_3)=0$;
\item[(4)] Reduction Lemma: For arbitrary subspace $K$ of $L_1\cap L_2+L_2\cap L_3+L_3\cap L_1$, $$s(L_1,L_2,L_3)=s(L_1^K,L_2^K,L_3^K),$$ where $L_i^K$ is the image of $L_i$ under the linear symplectic reduction $(K+K^\omega)\rightarrow {E}^K:=(K+K^\omega)/(K\cap K^\omega)$, where $K^\omega$ is the symplectic orthogonal complement of $K$ in $\mathbf{R}^{2n}$;
\item[(5)] The signature runs through all integers between $-\frac{1}{2}\dim {E}^F$ and $+\frac{1}{2}\dim {E}^F$, where $F=(L_1\cap L_2)+(L_2\cap L_3)+(L_3\cap L_1)$. Consequently, $s(L_1,L_2,L_3)+\dim(L_1\cap L_2)+\dim(L_2\cap L_3)+\dim(L_3\cap L_1)+n$ is even;
\item[(6)] The orbit of the action $\mathrm{Sp}(2n)$ on $\Lambda(n)^3$ are completely determined by $\dim(L_1\cap L_2\cap L_3),\dim(L_1\cap L_2),\dim(L_2\cap L_3),\dim(L_3\cap L_1),s(L_1,L_2,L_3)$: If these five numbers are the same for two triplets, then they lie on the same orbit;
\item[(7)] The signature is locally constant on the set of all triplets with given dimensions of intersections;
\item[(8)] symplectic additivity: for two symplectic vector space $(E_i,\omega_i)$ with corresponding triplets $(L_1,L_2,L_3),(\tilde L_1,\tilde L_2,\tilde L_3)$ respectively, we have $s(L_1\oplus \tilde L_1, L_2\oplus \tilde L_2, L_3\oplus \tilde L_3)=s(L_1,L_2,L_3)+s(\tilde L_1, \tilde L_2, \tilde L_3)$.
\end{itemize}
\end{proposition}

\subsection{Index for Pairs}

Given two continuous paths $L_1,L_2:[a,b]\rightarrow \Lambda(n)$ with $a\leq b$, we can choose a suitable partition $a=t_0\leq\cdots\leq t_k=b$ and Lagrangian subspaces $M_i$ such that $M_i$ is transversal to $L_1(t)$ and $L_2(t)$ for any $t_{i-1}\leq t\leq t_i$($i=1,...,k$), and define the {\bf intersection number} $[L_1:L_2]$ to be $$[L_1:L_2]:=\frac{1}{2}\sum_{i=1}^k \left(\tau(L_1(t_{i-1}),L_2(t_{i-1}),M_i)-s(L_1(t_i),L_2(t_i),M_i)\right).$$ One can prove that this definition is independent of the choices via cocycle identities and locally constancy of the signature.

\begin{proposition}\label{index for pairs}
The intersection number has the following properties:
\begin{itemize}
\item[(1)] Antisymmetry: $[L_1:L_2]+[L_2:L_1]=0$;
\item[(2)] Invariance: $[AL_1:AL_2]=[L_1:L_2]$ for any continuous path $A: [a,b]\rightarrow \mathrm{Sp}(2n)$;
\item[(3)] $[L_1:L_2]+\frac{1}{2}\dim(L_1(a)\cap L_2(a))+\frac{1}{2}\dim(L_1(b)\cap L_2(b))\in \mathbf{Z}$. In particular, $[L_1:L_2]$ is an integer if the intersections at the boundary are transversal;
\item[(4)] Given a third path $L_3: [a,b]\rightarrow \Lambda(n)$, $$[L_1:L_2]+[L_2:L_3]+[L_3:L_1]=\frac{1}{2}(s(L_1(a),L_2(a),L_3(a))-s(L_1(b),L_2(b),L_3(b)));$$
\item[(5)] The intersection number characterizes the connected components of the space of paths $L_1\times L_2:[a,b]\rightarrow \Lambda(n)^2$ with given intersection dimensions at the boundaries;
\item[(6)] If $K(t)$ is a continuous curve of isotropic subspaces contained in $L_1$ such that $\dim(K\cap L_2)$ is constant, then the reductions $L_i^K$ of $L_i$ with respect to $K$ are continuous, and $$[L_1:L_2]=[L_1^K:L_2^K].$$
\end{itemize}
\end{proposition}



\subsection{Leray Index}
\cite{Le81}, \cite{dG90}, \cite{CLM94}, \cite{dG09}

let $\pi: \tilde\Lambda(n)\rightarrow\Lambda(n)$ be the universal covering of the Lagrangian Grassmannian $\Lambda(n)$. Given $u_1,u_2\in \tilde\Lambda(n)$, choose a path $u:[a,b]\rightarrow \tilde\Lambda(n)$ connecting $u_1$ and $u_2$, and let $L(t)=\pi(u(t))$, then the {\bf Leray index} is defined to be $$\nu(u_1,u_2):=[L(t):L(b)]=[L(t):L(a)]\in\frac{1}{2}\mathbf{Z}.$$ One can prove that the definition is independent of the choice of the path $L(t)$ via Proposition \ref{index for pairs}.

\begin{proposition}\label{Leray index}
The Leray index has the following properties:
\begin{itemize}
\item[(1)] Leray formula: for $L_i=\pi(u_i), (i=1,2,3),$ $$\nu(u_1,u_2)+\nu(u_2,u_3)+\nu(u_3,u_1)=\frac{1}{2}s(L_1,L_2,L_3);$$
\item[(2)] For arbitrary lifts $u_i(t)$ of Lagrangian curves $L_i(t)$, $$[L_1:L_2]=\nu(u_1(a),u_2(a))-\nu(u_1(b),u_2(b));$$
\item[(3)] $\nu(u_1,u_2)$ is locally constant on the set of all pairs $(u_1,u_2)$ with fixed $\dim(L_1\cap L_2)$.
\end{itemize}
\end{proposition}

de Gosson (\cite{dG90}) proves that properties (1) and (3) imply that this definition of Leray index ia equivalent to the constructions in \cite{GS77} and \cite{LV80}.

\section{Maslov Index: Symplectic Group Perspective}

\subsection{Symplectic Groups}

We recall some topological facts on the symplectic group which is used in the definition of various Maslov indices from \cite{Long02} where the proofs and further information can be found.

For any $\omega$ in  the unit circle $\mathbf{U}$ in the complex plane, we denote by $\mathrm{Sp}(2n)_\omega^*$ the subset of $\mathrm{Sp}(2n)$ which contains all symplectic matrices without eigenvalue equal to $\omega$, and $\mathrm{Sp}(2n)_\omega^0$ its complement in $\mathrm{Sp}(2n)$. We call them the regular and the singular subsets of the symplectic group respectively, and we have $\mathrm{Sp}(2n)=\mathrm{Sp}(2n)_\omega^*\cup\mathrm{Sp}(2n)_\omega^0$. We can further decompose the singular subset into the disjoint union of $\mathcal{M}_\omega^k(2n)$ which consists of symplectic matrices with the complex dimension of the kernel of the matrix $M-\omega I_{2n}$ equal to $k$, and we have $\mathrm{Sp}(2n)_\omega^0=\cup_{k=1}^{2n} \mathcal{M}_\omega^k(2n)$. $\mathcal{M}_\omega^1(2n)$ is called the regular part of the singular subset $\mathrm{Sp}(2n)_\omega^0$.

The $\diamond$-product $M_1\diamond M_2$ of two real symplectic matrices of the square block form $$M_1=\left(\begin{matrix} A_1 & B_1 \cr C_1 & D_1\cr\end{matrix}\right)_{2i\times 2i},\,\, M_2=\left(\begin{matrix}A_2 & B_2 \cr C_2 & D_2\cr\end{matrix}\right)_{2j\times 2j}$$ is a $2(i+j)\times 2(i+j)$ symplectic matrix $$M_1\diamond M_2=\left(\begin{matrix}A_1 & 0 & B_1 & 0 \cr
                                      0   & A_2 & 0& B_2\cr
                                      C_1 &0 & D_1& 0 \cr
                                      0 & C_2 & 0 & D_2\cr\end{matrix}\right).$$ It is introduced to fit with the symplectic orthogonal decomposition of the standard symplectic vector space with respect to the standard symplectic matrix $J$. This product is associative, and closed in the symplectic group.

For the global structure of the symplectic group, we have the following

\begin{theorem}(Gelfand-Lidskii(\cite{GL55}); Moser(\cite{Moser58}))
The symplectic group $\mathrm{Sp}(2n)$ is homeomorphic to the topological product of of the unit circle $\mathbf{U}$ on the complex plane and a simply connected topological space.
\end{theorem}

The proof is based on the polar decomposition of symplectic matrices to get that the symplectic group $\mathrm{Sp}(2n)$ is homeomorphic to the product of the set of positive definite symmetric symplectic matrices and the orthogonal symplectic group. The former is homeomorphic to the contractible Euclidean space $\mathbf{R}^{n(n+1)}$, and the latter is isomorphic to the unitary group. It is well known that the unitary group is the product of the unit circle and the special unitary group which is simply connected following H. Weyl. As a corollary, we have

\begin{proposition} The fundamental group of the symplectic group is
 $\pi_1(\mathrm{Sp}(2n))\cong\mathbf{Z}$.
\end{proposition}

\begin{theorem}(Conley-Zehnder, Salamon-Zehnder, Long; c.f. \cite{Long02}, \S 2.4)\label{SpI} For any $\omega\in\mathbf{U}$,
the regular subset $\mathrm{Sp}(2n)_\omega^*$ possesses precisely two path connected components $\mathrm{Sp}(2n)_\omega^+:=\{M\in \mathrm{Sp}(2n)\,|\, (-1)^{n-1}\omega^{-n}\det(M-\omega I_{2n})<0\}$ and $\mathrm{Sp}(2n)_\omega^-:=\{M\in \mathrm{Sp}(2n)\,|\, (-1)^{n-1}\omega^{-n}\det(M-\omega I_{2n})>0\}$ both of which are simply connected in $\mathrm{Sp}(2n)$.
\end{theorem}

One can see easily that $$M_n^+:=\left(\begin{matrix} 2 & 0 \cr 0 &  \frac{1}{2}\cr\end{matrix}\right)^{\diamond n}\in \mathrm{Sp}(2n)_\omega^+$$ and $$M_n^-:=\left(\begin{matrix} -2 & 0 \cr 0 &  -\frac{1}{2}\cr\end{matrix}\right)\diamond\left(\begin{matrix} 2 & 0 \cr 0 &  \frac{1}{2}\cr\end{matrix}\right)^{\diamond (n-1)}\in \mathrm{Sp}(2n)_\omega^-$$ which will be used later on. In the following, for brevity, we use $D(\lambda)$ to denote the symplectic matrix of the form $$D(\lambda):=\left(\begin{matrix} \lambda & 0 \cr 0 &  \frac{1}{\lambda}\cr\end{matrix}\right),\,\,\lambda\neq 0.$$

\begin{theorem}(Long, c.f. \cite{Long02}, \S 2.5)\label{SpII} For any $n\geq 2$ and $\omega\in\mathbf{U}$, the singular subset $\mathrm{Sp}(2n)_\omega^0$ is path connected and not simply connected in $\mathrm{Sp}(2n)$ whence not simply connected in itself. $\mathrm{Sp}(2)_1^0$ and $\mathrm{Sp}(2)_{-1}^0$ are path connected and simply connected,  however $\mathrm{Sp}(2)_\omega^0$ has precisely two path connected components for $\omega\in\mathbf{U}\backslash\mathbf{R}$ each of which is simply connected.
\end{theorem}

\begin{theorem}(Long, c.f. \cite{Long02}, \S 2.7)\label{SpIII} For any $\omega\in\mathbf{U}$, the regular part $\mathcal{M}_\omega^1(2n)$ of the singular hypersurface $\mathrm{Sp}(2n)_\omega^0$ is a smooth codimension $1$ open submanifold in $\mathrm{Sp}(2n)$ and possesses a natural orientation.
\end{theorem}

The proofs are based on the the notion of basic normal form for symplectic matrices to be introduced later in the section.

{\bf $\mathrm{Sp}(2)$ model}

Here, we follow Long (1991, c.f. \cite{Long02}, \S 2.1) to introduce a geometric representation of $\mathrm{Sp}(2)$ in terms of $\mathbf{R}^3$-cylindrical coordinates. Be aware that there is another model due to Gelfand-Lidskii (\cite{GL55}). Long's model is suggestive and tailored for the index theory developed in the next subsection.
For any matrix $M\in \mathrm{Sp}(2)$, by the symplectic polar decomposition of a symplectic matrix into the product of a positive definite symmetric symplectic matrix and a symplectic orthogonal matrix, $M$ can be written in the following manner
$$M=\left(\begin{matrix} r & z \cr z &  \frac{(1+z^2)}{r}\cr\end{matrix}\right)\left(\begin{matrix}\cos\theta &-\sin\theta \cr \sin\theta & \cos\theta\cr\end{matrix}\right),$$ where $(r,\theta,z)\in \mathbf{R}^+\times (\mathbf{R}/2\pi\mathbf{Z})\times \mathbf{R}$ is uniquely determined by $M$. In fact, the map $\Phi: M\mapsto (r,\theta,z)$ defines a smooth diffeomorphism from $\mathrm{Sp}(2)$ to $\mathbf{R}^3\backslash\{z-\textrm{axis}\}$. Now for $\omega=\sin\varphi+\sqrt{-1}\cos\varphi$, $\mathrm{Sp}(2)_\omega^*=\mathrm{Sp}(2)_\omega^+\cup \mathrm{Sp}(2)_\omega^-$ with $$\mathrm{Sp}(2)_\omega^\pm=\{(r,\theta,z)\in \mathbf{R}^+\times (\mathbf{R}/2\pi\mathbf{Z})\times \mathbf{R}\,|\, \pm (r^2+z^2+1)\cos\theta>2r\cos\varphi\},$$ and
$$\mathrm{Sp}(2)_\omega^0=\{(r,\theta,z)\in \mathbf{R}^+\times (\mathbf{R}/2\pi\mathbf{Z})\times \mathbf{R}\,|\, (r^2+z^2+1)\cos\theta = 2r\cos\varphi\}$$ which can be further decomposed. Set $\mathrm{Sp}(2)_{\omega,\pm}^0=\{(r,\theta,z)\in \mathrm{Sp}(2)_\omega^0\,|\, \pm \sin\theta>0\}$. We have
$$\mathrm{Sp}(2)_{\pm 1}^0=\mathrm{Sp}(2)_{\pm 1,+}^0\cup\{\pm I_2\}\cup \mathrm{Sp}(2)_{\pm 1,-}^0,$$ and $$\mathrm{Sp}(2)_\omega^0=\mathrm{Sp}(2)_{\omega,+}^0\cup \mathrm{Sp}(2)_{\omega, -}^0.$$ All the above properties about $\mathrm{Sp}(2)$ can be red off directly from the beautiful  figures in \cite{Long02}(p 49).

{\bf Krein type}

Motivated by the linear stability problem in Hamiltonian systems with periodic coefficients, M. Krein (\cite{K50}, \cite{K51a}, \cite{K51b}, \cite{K51c}, \cite{K55}, see also \cite{YaS75}, \cite{E90} and \cite{Long02}) developed his stability theory of symplectic matrices and linear Hamiltonian systems which is also rediscovered independently by Moser (\cite{Moser58}), and Krein's theory plays a very important role in the development of the theory of Hamiltonian dynamics. We introduce the matrix $G=\sqrt{-1}J$ and use it to define the Krein quadratic form $\langle Gx,y\rangle,\,\forall x,y\in \mathbf{C}^{2n}$. If we denote by $$E_\lambda(M):=\cup_{k\geq 1}\ker_\mathbf{C}(M-\lambda I)^k\subset \mathbf{C}^{2n}$$ the complex root vector space/eigenspace, one can prove the following fact

\begin{proposition}
Suppose two eigenvalues $\lambda,\mu$ of a symplectic matrix is such that $\lambda\bar\mu\neq 1$, then $E_\lambda(M)$ and $E_\mu(M)$ are $G$-orthogonal.
\end{proposition}
So for distinct $\lambda,\mu\in\sigma(M)\cap\mathbf{U}$, their eigenspace must be mutually $G$-orthogonal. Furthermore, the restriction of the Krein form $G$ to $E_\lambda(M)$ is nondegenerate whence we can consider its total multiplicities of positive and negative eigenvalues which we denote by the pair $(p,q)$ and call it the {\bf Krein type number} of the eigenvalue $\lambda$. They will play a key role in the definition and computations of Conley-Zehnder index theory. If $q=0$, we say the eigenvalue is Krein positive. If $p=0$, we say it is Krein-negative. In both case, we say it is Krein definitive, otherwise mixed kind. To give an idea one can check directly that the matrix $R(\theta)=\left(\begin{matrix}\cos\theta &-\sin\theta \cr \sin\theta & \cos\theta\cr\end{matrix}\right)$ with $\theta\in(0,\pi)\cup(\pi,2\pi)$ has eigenvalues $\exp(\sqrt{-1}\theta)$ and $\exp(-\sqrt{-1}\theta)$ with Krein type $(0,1)$ and $(1,0)$ respectively. Some properties for further references are listed here. The Krein type numbers are conjugate invariant. If $\lambda$ has Krein type $(p,q)$, then $\bar\lambda$ has switched Krein type $(q,p)$. In particular, if $1$ or $-1$ is an eigenvalue, then its Krein type must be of form $(p,p)$ for some $p\in\mathbf{N}$.

\subsection{Conley-Zehnder/Maslov-type Index}

For $\tau >0$, we define the set of symplectic matrix path as $$\mathcal{P}_\tau(2n):=\{\gamma\in C([0,\tau],\mathrm{Sp}(2n))\,|\, \gamma(0)=I_{2n}\},$$
and the fundamental solution to a general linear Hamiltonian system with continuous symmetric periodic coefficients, which is the case for the linearization of a general Hamiltonian system along a periodic orbit, is an element $\gamma$ of $\mathcal{P}_\tau(2n)$. Since $\mathrm{Sp}(2n)$ is homeomorphic of a circle, $\gamma$ rotates naturally in $\mathrm{Sp}(2n)$ along the circle, and the Conley-Zehnder index encodes nicely the counting of rotations. Furthermore it fits very well with the variational nature of the orbit which we will see in \S \ref{Maslov=Morse}, and it can be seen as a generalization of the usual Morse index for closed geodesics on a Riemannian manifold.

Conley and Zehnder (\cite{CZ84}) defined the index theory named after them for nondegenerate paths in $\mathrm{Sp}(2n)$ with $n\geq 2$, and it is further defined by Long and Zehnder (\cite{LZ90}) for the nondegenerate paths in $\mathrm{Sp}(2)$. This index theory is extended to degenerate paths by Long and Viterbo independently (\cite{Long90},\cite{V90},\cite{Long97}). We shall use Conley-Zehnder index or Maslov-type index interchangeably in the following. Our main reference is the monograph by Y. Long (\cite{Long02}), and for a quick review with various applications please refer to \cite{Long08}.

Now we define the $\omega$-nullity and $\omega$-index for symplectic paths.

\begin{definition} For any $M\in \mathrm{Sp}(2n)$ and $\omega\in \mathbf{U}$, the {\bf $\omega$-nullity} $\nu_\omega(M)$ is defined to be $$\nu_\omega(M):=\dim_\mathbf{C}\ker_\mathbf{C}(M-\omega I_{2n}).$$ For any $\tau>0, \gamma\in \mathcal{P}_\tau(2n)$, we define the {\bf $\omega$-nullity} of a symplectic path $\gamma$ to be $$\nu_\omega(\gamma):=\nu_\omega(\gamma(\tau))=\dim_\mathbf{C}\ker_\mathbf{C}(\gamma(\tau)-\omega I_{2n}).$$ A path $\gamma\in\mathcal{P}_\tau(2n)$ is called {\bf $\omega$-degenerate} if $\nu_\omega(\gamma)>0$, otherwise it is called {\bf $\omega$-nondegenerate} the set of which is denoted by $\mathcal{P}_\tau^*(2n)$.
\end{definition}

For simplicity, we firstly give the definition of Maslov-type index $i_\omega(\gamma)$ in the case of $\omega=1$.

\begin{definition}
Given two paths $\gamma_0$ and $\gamma_1\in \mathcal{P}_\tau(2n)$, if there is a map
$\delta\in C([0, 1] \times [0, ¦Ó ], \mathrm{Sp}(2n))$ such that $\delta(0, \cdot) =\gamma_0(\cdot), \delta(1, \cdot) =\gamma_1(\cdot), \delta(s, 0) = I_{2n}$, and $\nu_\tau(\delta(s, \cdot))$ is
constant for $0 \leq s \leq 1$, then we say that $\gamma_0$ and $\gamma_1$ are homotopic on $[0, \tau ]$ along $\delta(\cdot, \tau )$ and we write $\gamma_0\sim\gamma_1$
on $[0, \tau ]$ along $\delta(\cdot, \tau )$. This homotopy possesses fixed end points if $\delta(s, \tau ) =\gamma_0(\tau )$ for all $s \in [0, 1]$.
\end{definition}

As recalled above, any symplectic matrix $M$ has the unique polar decomposition $M = AU$, where $A =
(MM^T )^{1/2}$ is symmetric positive definite and symplectic, $U$ is orthogonal and symplectic.
Now $U$ must have the form
$$U =\left(\begin{matrix} u_1 & -u_2 \cr u_2 & u_1\cr\end{matrix}\right),$$
where $u = u_1+\sqrt{-1}u_2$ is unitary. So for every path $\gamma\in \mathcal{P}_\tau(2n)$, we can associate
a path $u(t)$ in the unitary group on $\mathbf{C}_n$ to it. If $\Delta(t)$ is any continuous real function satisfying $\det u(t) = \exp(\sqrt{-1}\Delta(t))$, the difference $\Delta(\tau)-\Delta(0)$ depends only on $\gamma$ but not on the choice of
the function $\Delta(t)$. Therefore we may define the mean rotation number of $\gamma$ on $[0, \tau ]$ by
$$\Delta_\tau(\gamma) = \Delta(\tau )-\Delta(0).$$

\begin{proposition}
 If $\gamma_0$ and $\gamma_1\in\mathcal{P}_\tau(2n)$ possess common end point $\gamma_0(\tau) =
\gamma_1(\tau )$, then $\Delta_\tau (\gamma_0) =\Delta_\tau (\gamma_1)$ if and only if $\gamma_0\sim\gamma_1$ on $[0, \tau]$ with fixed end points.
\end{proposition}

For any nondegenerate path $\gamma\in \mathcal{P}_\tau^*(2n)$, there exists a path $\beta : [0, \tau ]\rightarrow \mathrm{Sp}^*(2n)$ such that
$\beta(0) =\gamma(\tau)$ and $\beta(\tau) = M_n^+$ or $M^-_n$. Define their concatenation product path $\beta*\gamma$ by
$$\beta*\gamma(t) =\left\{\begin{matrix}\gamma(2t),\,\,\,\, 0\leq t\leq \tau/2 \cr
                                         \beta(2t-\tau),\,\,\,\, \tau/2\leq t\leq \tau\end{matrix}\right. $$
Then $k\equiv \Delta_\tau(\beta*\gamma)/\pi\in \mathbf{Z}$ and is independent of the choice of the path $\beta$ by the simply connectedness of each component of $\mathrm{Sp}(2n)^*$. In this case we write $\gamma\in \mathcal{P}_{\tau,k}^*(2n)$.

\begin{proposition}
$\mathcal{P}_{\tau,k}^*(2n)$'s give a homotopy classification of $\mathcal{P}_\tau^*(2n)$.
\end{proposition}

\begin{definition}If $\gamma\in\mathcal{P}_{\tau,k}^*(2n)$, we define $i_1(\gamma) = k$.
\end{definition}
In fact, in each homotopy class one can define some standard zigzag nondegenerate symplectic path using only very simple $2\times 2$ symplectic matrices and $\diamond$-product(\cite{Long02}, p112).

The following proposition is fundamental in the definition of the Maslov-type index for degenerate symplectic paths.

\begin{proposition}
For any ¦Ã$\gamma\in\mathcal{P}^0_\tau(2n)$, there exists a one parameter
family of symplectic paths $\gamma_s$ with $s\in [-1, 1]$ and a $t_0 \in (0, \tau)$ sufficiently close to $\tau$ such that
\begin{itemize}
\item[(i)] $\gamma_0=\gamma,\,\,\gamma_s(t)=\gamma(t)\,\,\textrm{for}\,\,0\leq t\leq t_0$;
\item[(ii)] $\gamma_s\in \mathcal{P}_\tau^*(2n)$\,\,\textrm{for}\,\, $s\in[-1,1]\backslash\{0\}$;
\item[(iii)] $i_1(\gamma_s)=i_1(\gamma_{s'})$,\,\,\textrm{if}\,\, $ss'>0$;
\item[(iv)] $i_1(\gamma_1)-i_1(\gamma_{-1}=\nu_1(\gamma)$;
\item[(v)] $\gamma_s\rightarrow \gamma_0=\gamma$\,\,\textrm{in}\,\,$\mathcal{P}_\tau(2n)$\,\,\textrm{as}\,\,$s\rightarrow 0$.
\end{itemize}
\end{proposition}

\begin{definition}
$i_1(\gamma):=i_1(\gamma_s)$ for $s\in [-1,0)$.
\end{definition}

In fact, we have the following theorem to show that the definition is independent of the choice.

\begin{theorem}
For any $\gamma\in\mathcal{P}_\tau^0(2n)$, and any $\beta\in\mathcal{P}_\tau^*(2n)$ which
is sufficiently close to $\gamma$, there holds
$$i_1(\gamma)=i_1(\gamma_{-1})\leq i_1(\beta)\leq i_1(\gamma_1)=i_1(\gamma)+\nu_1(\gamma).$$
\end{theorem}

This theorem gives another way to characterize the Maslov-type index for degenerate path. Namely, $$i_1(\gamma)=\inf\{i_1(\beta)\,|\, \beta\in\mathcal{P}_\tau^*(2n),\,\,\beta\,\, is\,\, sufficiently\,\, close\,\, to\,\, \gamma\,\, in\,\, \mathcal{P}_\tau(2n)\}.$$

So for any path $\gamma\in\mathcal{P}_\tau(2n)$, we have defined a pair of integers $(i_1(\gamma),\nu_1(\gamma))\in \mathbf{Z}\times\{0,1,...,2n\}$ which we call the {\bf Maslov-type index} of the symplectic path.

We have the following axiomatic characterization for Maslov-type index on any continuous path in $\mathcal{P}_\tau(2n)$.

\begin{theorem}The Maslov-type index $i_1: \cup_{n\in \mathbf{N}}\mathcal{P}_\tau(2n)\rightarrow \mathbf{Z}$ is uniquely determined by the following five axioms:
\begin{itemize}
\item[(i)]({\bf Homotopy Invariance}) For $\gamma_0$ and $\gamma_1\in \mathcal{P}_\tau(2n)$ such that $\gamma_0\sim\gamma_1$ on $[0,\tau]$, $$i_1(\gamma_0)=i_1(\gamma_1);$$
\item[(ii)] ({\bf Symplectic Additivity}) For any $\gamma_i\in\mathcal{P}_\tau(2n_i)(i = 0,1)$, there holds
               $$i_1(\gamma_0\diamond\gamma_1)=i_1(\gamma_0)+i_1(\gamma_1);$$
\item[(iii)]({\bf Clockwise Continuity}) For any $\gamma\in \mathcal{P}_\tau^0(2)$ with $\gamma(\tau)=N_1(1,b)=\left(\begin{matrix} 1 & b \cr 0 &  1\cr\end{matrix}\right)$ with $b=\pm 1$ or $0$, there exists a $\theta_0>0$ such that $$i_1([\gamma(\tau)R(-\theta \frac{t}{\tau})]*\gamma)=i_1(\gamma),\,\,\,\, \forall 0<\theta\leq\theta_0; $$
\item[(iv)] ({\bf Counterclockwise Jumping}) For any $\gamma\in \mathcal{P}_\tau^0(2)$ with $\gamma(\tau)=N_1(1,b)$ with $b=\pm 1$, there exists a $\theta_0>0$ such that $$i_1([\gamma(\tau)R(\theta \frac{t}{\tau})]*\gamma)=i_1(\gamma)+1,\,\,\,\, \forall 0<\theta\leq\theta_0; $$
\item[(v)] ({\bf Normality}) For the standard path $\alpha=D(1+\frac{t}{\tau}), t\in[0,\tau]$, there holds $$i_1(\alpha)=0.$$
\end{itemize}
\end{theorem}


Based on Theorems \ref{SpI}, \ref{SpII} and \ref{SpIII}, one can define $\omega$-Maslov-type index using similar elementary methods above. Especially, for degenerate path $\gamma\in\mathcal{P}_{\tau,\omega}^0(2n)$, the $\omega$-index is defined to be
$$i_\omega(\gamma):=\inf \{i_\omega(\beta)\,|\,\beta\in \mathcal{P}_{\tau,\omega}^*(2n)\,\, and\,\, \beta\,\, is\,\, sufficiently\,\, close\,\, to\,\, \gamma\,\, in\,\, \mathcal{P}_\tau(2n)\}.$$ When $\omega=1$, the $\omega$-index theory coincides with the Maslov-type index theory. An axiom characterization of the $\omega$-index theory can also be given (\cite{Long02},p147).

As a corollary, we can give an axiomatic characterization of the Maslov-type index for any end point free curve in the symplectic group. For any curve $f\in C([a,b],\mathrm{Sp}(2n))$, choose $\xi\in \mathcal{P}_1(2n)$ so that $\xi(1)=f(a)$. Define the new path $\eta$ to be the concatenation of $\xi$ and $f$. Now both $\xi$ and $\eta$ are paths in $\mathcal{P}_1(2n)$. Then we define
   $$i(f):=i_1(\eta)-i_1(\xi).$$ It is easy to see that the definition only depends on $f$ itself and is well defined. Note that here we do not require the end points of the path $f$ to be nondegenerate. Now we have the axiomatic characterization for this index

\begin{theorem}(\cite{Long02}, p148) The index $i$ defined above for continuous curves in the symplectic group is uniquely determined by the following axioms:
\begin{itemize}
\item[(i)] ({\bf Homotopy Invariance}) Two continuous curves in $\mathrm{Sp}(2n)$ with the same initial and end points possess the same index if and only if they can be continuously deformed to each other with the end points fixed;
\item[(ii)] ({\bf Vanishing}) $i(f)=0$ if $f\in C([a,b],\mathrm{Sp}(2n))$ with $\nu_1(f(t))=$constant;
\item[(iii)] ({\bf Symplectic Additivity}) $i(f_0\diamond f_1)=i(f_0)+i(f_1)$ if $f_i\in C([a,b],\mathrm{Sp}(2n_i))$(i=0,1);
\item[(iv)] ({\bf Concatenation}) $i(f)=i(f|_{[a,b]})+i(f|_{[b,c]})$ if $f\in C([a,c],\mathrm{Sp}(2n))$ such that $a<b<c$;
\item[(v)] ({\bf Normalization}) For $f(t)=N_1(1,b)R(t/2)$ with $t\in [-1,1]$ and $b=\pm 1$ or $0$, $$i(f|_{[-1,0]})=0,\,\,\,\,i(f|_{[0,1]})=2-|b|.$$
\end{itemize}
\end{theorem}

Motivated by the minimality of period, multiplicity and stability problems of geometrically distinct closed geodesics in Riemannian geometry, R. Bott (\cite{Bott56}) initiated the studies on the iteration formula for Morse index theory in 1956. It is a natural question to consider the iteration formula for Maslov index theory of closed characteristics of Hamiltonian systems to which many works are devoted. In fact, it is a symbiosis with the establishment of the Maslov-type index theory. And we refer the reader to \cite{Long02}(p. xiii) for the detailed history.

For any path $\gamma\in\mathcal{P}_\tau(2n)$ and $m\in\mathbf{N}$, we denote the $m$-th iteration $\gamma^m\in\mathcal{P}_{m\tau}(2n)$ of $\gamma$ by
$$\gamma^m(t)=\gamma(t-j\tau)\gamma(\tau)^j,\,\, j\tau\leq t\leq (j+1)\tau,\,\, j=0,1,...,m-1,$$ and we denote the index and nullity of $\gamma^m$ by $$(i(\gamma,m),\nu(\gamma,m))=(i_1(\gamma^m),\nu_1(\gamma^m)).$$

\begin{theorem}(Bott-type formula, c.f. \cite{Long02}, Th.9.2.1) For any $\tau>0$, $\gamma\in\mathcal{P}_\tau(2n)$, $z\in \mathbf{U}$ and $m\in \mathbf{N}$,
\bes i_z(\gamma,m)=\sum_{\omega^m=z}i_\omega(\gamma),\\
     \nu_z(\gamma,m)=\sum_{\omega^m=z}\nu_\omega(\gamma).
\ees
\end{theorem}

To get a rough idea of the proof, we first recall the inverse homotopy theorem which is very useful in the study of the iteration theory for Maslov-type indices.

\begin{theorem}
For any two paths $\gamma_0$ and $\gamma_1\in \mathcal{P}_\tau(2n)$ with $i_1(\gamma_0) = i_1(\gamma_1)$, suppose that there exists a continuous path $h : [0, 1]\rightarrow \mathrm{Sp}(2n)$ such that $h(0) = \gamma_0(\tau )$, $h(1) = \gamma_1(\tau)$, and $\dim\ker(h(s)-I) = \nu_\tau(\gamma_0)$ for all $s \in [0, 1]$. Then $\gamma_0\sim\gamma_1$ on $[0, \tau]$ along $h$.
\end{theorem}

The proof of the Bott-type iteration formula goes as follows. For the fundamental solution of the linearized Hamiltonian system along a periodic solution, certain corresponding differential operator can be made in diagonal form with respect to the orthogonal decomposition related to the iteration, then the relation between the (relative) Morse index and Maslov-type index to be indicated in the next section is resorted to get the result.

Such a formula works for the fundamental solutions of general linear Hamiltonian systems with continuous symmetric and periodic coefficients, and hence for periodic solutions of any nonlinear Hamiltonian systems with periodic in time Hamiltonian functions.

Various generalizations exist in the literature. For example, finite group theoretic point of view is given in \cite{HS09} where we understand the iteration as a finite cyclic group action on the closed paths and extended the Bott-type iteration formula to more general finite groups.

For a given path $\gamma$ of symplectic matrices originating from some periodic orbits, in general it is very difficult to compute the above Maslov indices by definitions and effective tools for computations are desirable. Y. Long successfully built up a method for this purpose. Here the goal is to replace $\gamma$ by a new path $\eta\in \mathcal{P}_\tau(2n)$ such that
\bea i_1(\gamma^m)=i_1(\eta^m),\,\,\,\,\nu_1(\gamma^m)=\nu_1(\eta^m),\,\,\,\,\forall m\in \mathbf{N}
 \label{indexequality}\eea
and $i_1(\eta^m)$ and $\nu_1(\eta^m)$ are easier to compute. The natural idea is to construct a homotopy $\delta_s(t)(s\in [0,1],t\in[0,\tau])$ between $\gamma$ and $\eta$ in $\mathcal{P}_\tau(2n)$ with the end points $\delta_s(\tau)$ in certain subset $\Gamma$ of $\mathrm{Sp}(2n)$ with constraints (\ref{indexequality}) above. The traditional choice for $\Gamma$ is the conjugacy class of $\gamma(\tau)$ in $\mathrm{Sp}(2n)$ whence one gets new paths with end point matrices the symplectic normal forms of $\gamma(\tau)$. As is known (see, e.g., \cite{Long02}, Chapter 1) that the order of such normal forms can be as large as that of $\gamma(\tau)$, and the computations are still complicated. So a larger $\Gamma$ is needed. (\ref{indexequality}) is equivalent that the index function $i_\omega(\delta_s(\cdot)$ unchange for all $s\in[0,1]$ whenever $\omega$ is the spectrum of $\gamma(\tau)$ on the unit circle $\mathbf{U}$ by Bott-type iteration formula, which in turn is the same as the fact that the nullity function $\nu_\omega(\delta_s(\tau))$ of $s$ is constant whenever $\omega$ is the eigenvalue of $\gamma(\tau)$ on $\mathbf{U}$ by the stability property of Morse indices of symmetric matrices (\cite{Long02}, Lemma 6.1.3) and the fact that all the roots of unity are dense in $\mathbf{U}$. This motivated Long to introduce the following set containing the conjugate set of $\gamma(\tau)$ $$\Omega(\gamma(\tau)):=\{N\in \mathrm{Sp}(2n)\,|\, \sigma(N)\cap \mathbf{U}=\sigma(\gamma(\tau))\cap\mathbf{U},\,\,\nu_\lambda(N)=\nu_\lambda(\gamma(\tau)),\,\forall \lambda\in \sigma(\gamma(\tau))\cap\mathbf{U}\}$$ and the set $\Gamma$ is defined to be the path connected component $\Omega^0(\gamma(\tau))$ of $\Omega(\gamma(\tau))$ containing $\gamma(\tau)$ which is called the homotopy component of $\gamma(\tau)\in \mathrm{Sp}(2n)$. It is clear that $\Omega^0(\gamma(\tau))$ is the largest possible set for this purpose.

With the help of $\Omega^0(\gamma(\tau))$, we can decompose the symplectic normal forms further to much simpler {\bf basic normal form} (\cite{Long99}, \cite{Long02}(\S 1.8)). As we will see, they simplify a lot the computations of Maslov-type indices and unravel many significant features of symplectic group. More precisely, They are some special $2\times 2$ and $4\times 4$ symplectic matrices which we list below. The basic normal forms for eigenvalues outside of $\mathbf{U}$ are $$D(2)=\left(\begin{matrix}2 & 0 \cr 0 & \frac{1}{2}\cr\end{matrix}\right),\,\,\,\, D(-2)=\left(\begin{matrix} -2 & 0 \cr 0 & -\frac{1}{2}\cr\end{matrix}\right)\in \mathrm{Sp}(2).$$ The basic normal forms for eigenvalues $\omega=\exp(\sqrt{-1}\theta)\in\mathbf{U}$ are
\bes
N_1(\lambda,b)= \left(\begin{matrix}\lambda & b \cr 0 & \lambda\cr\end{matrix}\right) \in \mathrm{Sp}(2),\,\,\lambda=\pm 1,\,\,b=\pm 1,0;\\
R(\theta)=\left(\begin{matrix} \cos\theta & -\sin\theta \cr \sin\theta & \cos\theta\cr\end{matrix}\right)\in \mathrm{Sp}(2),\,\,\theta\in (0,\pi)\cup(\pi,2\pi);\\
N_2(\omega, b)=\left(\begin{matrix} R(\theta) & b \cr 0 & R(\theta)\cr\end{matrix}\right)\in \mathrm{Sp}(4),\,\, \theta\in (0,\pi)\cup(\pi,2\pi);\\
b=\left(\begin{matrix} b_1 & b_2 \cr b_3 & b_4\cr\end{matrix}\right),\,\, b_i\in \mathbf{R}, b_2\neq b_3.
\ees

\begin{theorem}(\cite{Long02}, Th. 2.3.8)\label{basicnormalformdecomposition} For any $M\in \mathrm{Sp}(2n)$, there is a path $f:[0,1]\rightarrow \Omega^0(M)$ such that $f(0)=M$ and $$f(1)=M_1(\omega_1)\diamond \cdots\diamond M_k(\omega_k)\diamond M_0,$$ where each $M_i(\omega_i)$ is a basic normal form of the eigenvalue $\omega_i\in \mathbf{U}$ for $1\leq i\leq k$, and the symplectic matrix $M_0=D(2)^{\diamond j}$ or $D(-2)\diamond D(2)^{\diamond (j-1)}$ for some nonnegative integer $j$.
\end{theorem}

With Bott-type iteration formula in mind, the properties of the index function $i_\omega(\gamma)$ of $\gamma$ at $\omega$ as a function of $\omega\in\mathbf{U}$ plays a key role. We follow \cite{Long02} (\S 9.1) to introduce the splitting number of a symplectic matrix at some $\omega\in\mathbf{U}$ which is used to detect the possible jumps of the index function at $\omega$.

For a fixed path $\gamma\in\mathcal{P}_\tau(2n)$, the index $i_\omega(\gamma)$ is a step function of $\omega\in\mathbf{U}$ which is constant on $\mathbf{U}\backslash \sigma(\gamma(\tau))$ with possible jumps only at eigenvalues of $\gamma(\tau)$ on $\mathbf{U}$. This motivates us to introduce the splitting numbers $S_M^\pm(\omega)$ to measure this jump.

\begin{definition} For any $M\in \mathrm{Sp}(2n)$, and $\omega\in\mathbf{U}$, choosing $\tau>0$ and $\gamma\in\mathcal{P}_\tau(2n)$ with $\gamma(\tau)=M$, we define the splitting numbers of $M$ at $\omega$ by
$$S_M^\pm(\omega):=\lim_{\epsilon\rightarrow 0^+}i_{\exp(\pm\epsilon\sqrt{-1})\omega}(\gamma)-i_\omega(\gamma).$$
\end{definition}

One can prove that the splitting numbers are well defined, in other words, they are independent of the choice of the symplectic path. Moreover, for $\omega\in\mathbf{U}$ and $M\in \mathrm{Sp}(2n)$, $S^\pm_N(\omega=S^\pm_M(\omega)$, for any symplectic matrix $N$ in the homotopy component $\Omega^0(M)$ of $M$.

We have the following axiomatic characterization of splitting numbers.

\begin{theorem}(\cite{Long02}, p198)
The integer valued splitting number pair $(S_M^+(\omega),S_M^-(\omega))$ defined for all $(\omega,M)\in \mathbf{U}\times\cup_{n\geq 1}\mathrm{Sp}(2n)$ are uniquely determined by the following axioms:
\begin{itemize}
\item[(1)] (Homotopy invariance) $S_M^\pm(\omega)=S_N^\pm(\omega)$ for all $N\in\Omega^0(M)$;
\item[(2)] (Symplectic additivity) $S_{M_1\diamond M_2}^\pm(\omega)=S_{M_1}^\pm(\omega)+S_{M_2}^\pm(\omega)$;
\item[(3)] (Vanishing) $S_{M}^\pm(\omega)=0$ if $\omega\notin\sigma(M)$;
\item[(4)] (Normalization) $(S_{M}^+(\omega),S_{M}^-(\omega))$ coincides with the ultimate type of $\omega$ for $M$ when $M$ is any basic normal form as listed above.
\end{itemize}
\end{theorem}

To get a better understanding on the splitting numbers, we need some further properties of the spectrum of the symplectic matrix, namely the {\bf ultimate type} of basic normal form which plays a fundamental role for the Maslov-type index and its iterations. A basic normal form matrix $M\in \mathrm{Sp}(2n)$ is called trivial if $MR((t-1)\alpha)^{\diamond n}$ possesses no eigenvalues on $\mathbf{U}$ for $t\in [0,1)$, and is non-trivial otherwise. Among all the basic normal forms, the matrices $N_1(1,-1), N_1(-1,1), D(2), D(-2), N_2(\omega,b)$ and $N_2(\bar\omega,b)$ with $\omega=\exp(\sqrt{-1}\theta)\in\mathbf{U}\backslash\mathbf{R}$ and $(b_2-b_3)\sin\theta>0$ are trivial, and the others are non-trivial. We define the ultimate type $(p,q)$ of $\omega\in \sigma(M)\cap \mathbf{U}$ for any basic normal form matrix $M$ to be its Krein type if $M$ is non-trivial, otherwise to be $(0,0)$ if $M$ is trivial. We also define the ultimate type of $\omega$ for $M$ to be $(0,0)$ if $\omega\in \mathbf{U}\backslash\sigma(M)$. Using the basic normal form decomposition of symplectic matrix (Theorem \ref{basicnormalformdecomposition}), one can finally define the {\bf ultimate type} of $\omega\in\mathbf{U}$ for $M\in \mathrm{Sp}(2n)$ to be the summation of the  ultimate types of $\omega$ for each basic normal form factor other that $M_0$. It is clear that the ultimate type is well defined, uniquely determined by $\omega$ and $M$, and is constant on the homotopy component of $M$ for fixed $\omega\in\mathbf{U}$. A basic observation about the Krein types is that whenever the eigenvalue $\omega$ leaves $\mathbf{U}$ by a small perturbation on a symplectic matrix $M$, both the Krein positive and negative type numbers of $\omega$ decrease by the same amount. Combined with the concrete computations for the basic normal forms, one can see that the differences between the Krein type and ultimate type for the positive part and the negative part are the same!

The interesting relation between the algebraically defined ultimate type and the splitting number is the following

\begin{theorem}(\cite{Long02}, p192) For any $\omega\in\mathbf{U}$ and $M\in \mathrm{Sp}(2n)$, $$S_M^+(\omega)=p,\,\,\,\, S_M^-(\omega)=q,$$ where $(p,q)$ is the ultimate type of $\omega$ for $M$.
\end{theorem}

A corollary of the theorem is that the index jump $\lim_{\epsilon\rightarrow 0^+}(i_{\exp(\epsilon\sqrt{-1}\omega)}(\gamma)-i_{\exp(-\epsilon\sqrt{-1}\omega)}(\gamma))$ at some $\omega\in\mathbf{U}\cap\sigma(\gamma(\tau))$ is the difference between the Krein positive and negative type number which is also the difference of the positive and negative ultimate number. Also one can see that the splitting numbers are nonnegative integers bounded above by the nullity $\nu_\omega(M)$ of $M$ at $\omega$. In fact, using the basic normal forms we have more precise estimates in terms of Krein type numbers $(p_\omega(M),q_\omega(M))$ $$0\leq \nu_\omega(M)-S_M^-(\omega)\leq p_\omega(M),\,\, 0\leq \nu_\omega(M)-S_M^+(\omega)\leq q_\omega(M).$$

Combining all these notions and theorems together, we can reduce the computations of splitting numbers of any general symplectic matrix to those of basic normal forms which are easier and listed below for completeness.

{\bf Splitting numbers of basic normal forms} (c.f. \cite{Long02}, p198)

\begin{itemize}
\item[(1)] $(S_M^+(1), S_M^-(1))=(1,1)$ for $M=N_1(1,b)$ with $b=1$ or $0$;
\item[(2)] $(S_M^+(1), S_M^-(1))=(0,0)$ for $M=N_1(1,-1)$;
\item[(3)] $(S_M^+(-1), S_M^-(-1))=(1,1)$ for $M=N_1(-1,b)$ with $b=-1$ or $0$;
\item[(4)] $(S_M^+(-1), S_M^-(-1))=(0,0)$ for $M=N_1(-1,1)$;
\item[(5)] $(S_M^+(\exp(\sqrt{-1}\theta)), S_M^-(\exp(\sqrt{-1}\theta)))=(0,1)$ for $M=R(\theta)$ with $\theta\in (0,\pi)\cup(\pi,2\pi)$;
\item[(6)] $(S_M^+(\omega), S_M^-(\omega))=(1,1)$ for $M=N_2(\omega,b)$ being nontrivial with $\omega=\exp(\sqrt{-1}\theta)\in \mathbf{U}\backslash\mathbf{R}$;
\item[(7)] $(S_M^+(\omega), S_M^-(\omega))=(0,0)$ for $M=N_2(\omega,b)$ being trivial with $\omega=\exp(\sqrt{-1}\theta)\in \mathbf{U}\backslash\mathbf{R}$;
\item[(8)] $(S_M^+(\omega), S_M^-(\omega))=(0,0)$ for any $\omega\in \mathbf{U}$ and $M\in \mathrm{Sp}(2n)$ such that $\sigma(M)\cap \mathbf{U}=\emptyset$.
\end{itemize}

Using the splitting numbers, we can state an abstract precise iteration formula.

\begin{theorem}(\cite{Long02}, \S 9.3) For any path $\gamma\in \mathcal{P}_\tau(2n)$ with $M=\gamma(\tau)$ and any natural number $m$, we have
\bes i(\gamma,m)&=& m(i(\gamma,1)+S_M^+(1)-C(M))\\
                &+& 2\sum_{0<\theta<2\pi}\left[\frac{m\theta}{2\pi}\right]S_M^-(\exp(\sqrt{-1}\theta))-(S_M^+(1)+C(M))
\ees
where $C(M)=\sum_{0<\theta<2\pi}S_M^-(\exp(\sqrt{-1}\theta))$ and $\left[\frac{m\theta}{2\pi}\right]$ denotes the least integer great than or equal to $\frac{m\theta}{2\pi}$.
\end{theorem}

Based on the Theorem \ref{basicnormalformdecomposition} about the basic normal form decomposition and their splitting numbers, we have a very concrete precise iteration formula (\cite{Long08}, Theorem 6.1, p24).  Based on this, one can get more precise information on mean index and sharper iteration inequalities (\cite{Long08} and references therein).

The merits of Maslov-type index theory develpoed by Y. Long and his collaborators are as follows. The usual Maslov index is defined only for loops or curves with nondegenerate endpoints which is enough for some topological issues, and must be extended to general paths allowing degenerate ends which is crucial for dynamical problems. The relations between Morse index and Maslov index are also extended from the usual Lagragian systems to more general nonlinear Hamiltonian systems which will be discussed in more details in \S \ref{Maslov=Morse}. Moreover, the iteration theory is completely established which is an essential ingredient in Gutzwiller's semiclassical trace formula.

\subsection{Relation Between Two Perspectives}

Recall the basic fact that the graph $$\Gamma_A:=\{(Ax,x)\,|\, x\in \mathbf{R}^{2n}\}$$ of a symplectic linear transformation $A\in \mathrm{Sp}(2n)$ is a Lagrangian subspace of $(\mathbf{R}^{2n}\times \mathbf{R}^{2n}, \textrm{pr}_1^*\omega-\textrm{pr}_2^*\omega)$. For a path $\gamma\in\mathcal{P}_\tau(2n)$, $i_1(\gamma)=[\Delta:\Gamma_{\gamma(t)}]$ with $\Delta$ the diagonal in $\mathbf{R}^{2n}\times \mathbf{R}^{2n}$ (\cite{CZ84}). More generally, we have, for any $L,M\in\Lambda(n)$,
\bea
i_1(\gamma)=[M: \gamma(t)L]+\frac{1}{2}s(\Delta, L\times M, \Gamma_{\gamma(\tau)}).\label{CZ=Maslovpahse}
\eea It follows from Proposition \ref{index for pairs} (4) by noting that $[M:\gamma(t)L]=[M\times L:\Gamma_{\gamma(t)}]$.

Cappell-Lee-Miller(\cite{CLM94}) considered the relations among various versions of Maslov index from geometry, topology and analysis viewpoints. It is still a fascinating topic to compare the variants of Maslov index bearing in mind important progresses in each definitions, especially in the field of Hamiltonian systems since the publication of \cite{CLM94} and the huge physical literature devoted to the computations. A good idea to show the equivalences and relations is by using axiomatic characterization of the indices as did in \cite{CLM94} and \cite{Long02}.

\section{Maslov Index and Morse Index}\label{Maslov=Morse}

\cite{Long02}, \cite{Ab01}

Maslov index for periodic orbits of a general Hamiltonian system is a kind of finite dimensional representation of infinite Morse index since the periodic orbit is the  critical point of the strongly indefinite action functional corresponding to the Hamiltonian system whose study was pioneered by Rabinowitz in 1978 (\cite{Ra78}). There are various definitions for Maslov index in the literature, and each of them is well adapted for the corresponding problems and has their own merits. The key point for Maslov-type index is its relation to Morse index theory for general nonlinear Hamiltonian systems. The theme here is that some relative Morse index is equal to Maslov-type index which can be proved, for example via spectral flow.

Maslov index is in phase space formulation (our main symplectic point in this paper), whereas Morse index is in configuration formulation which is more in the spirit of calculus of variations and physics.

Considering the periodic boundary value problem of the following Hamiltonian system  \bea \dot{x}(t)=JH'(t,x(t)),\,\,\,\, x(\tau)=x(0)\label{HS}\eea where $H\in C^2((\mathbf{R}/\tau\mathbf{Z})\times \mathbf{R}^{2n},\mathbf{R})$ for some fixed $\tau>0$ such that $\|H\|_{C^2}$ is finite and $H'(t,x)$ denotes the derivative with respect to the variable $x$. It is well known that the $\tau$-periodic solutions are in one-to-one correspondence with the critical points of the following action functional $$f(x)=\int_0^\tau (-\frac{1}{2}J\dot{x}\cdot x-H(t,x(t)))dt$$ for $x\in\textrm{dom}(A)\subset L_\tau\equiv L^2(\mathbf{R}/\tau\mathbf{Z}, \mathbf{R}^{2n})$ with the operator $A=-J\frac{d}{dt}$. The Morse indices of $f$ at its critical point are those of the following quadratic form on $L_\tau$
$$\phi(y):=\int_0^\tau (-J\dot{y}\cdot y-B(t)y\cdot y)dt$$ with $B(t)=H''(t,x(t))$ the symmetric matrix function of $t$ along the periodic solution. However these indices are infinite. One way to surround this difficulty is by using the saddle point reduction method on $L_\tau$ (\cite{AZ80},\cite{C93},\cite{Long02}). More precisely, one can get a finite dimensional subspace $Z\subset L_\tau$ consisting of finite Fourier polynomials with $\dim Z=2d$ large enough, an injective map $u:Z\rightarrow L_\tau$ and a function $z: Z\rightarrow \mathbf{R}$ such that there holds $$a(z)=f(u(z)),\,\,\,\,\forall z\in Z$$ and that the critical points of $a$ and $f$ are in one-to-one correspondence. Now we have the following important theorem.

\begin{theorem}(Conley-Zehnder\cite{CZ84}, Long-Zehnder\cite{LZ90}, Long\cite{Long02})
With the above notations, let $z$ be a critical point of the function $a$ and $x = u(z)$ be the corresponding solution of the Hamiltonian system (\ref{HS}). Denote the
Morse indices of $a$ at $z$ by $m^*(z)$ for $* = +, 0, -$. Then the Maslov-type index $(i_1(x), \nu_1(x))$ satisfy
$$m^-(z) = d + i_1(x), m^0(z) = \nu_1(x), m^+(z) = d - i_1(x) -\nu_1(x).$$
\end{theorem}

According to the theorem, the Maslov-type indices can be viewed as a finite representation of the
infinite Morse indices which captures the essential information of the variational problem. We should also note that for general Hamiltonian
 whose second derivative may not be bounded, similar theorem can be derived by the Galerkin approximations as did in \cite{FQ96}.

In configuration space formulation, one can consider similar periodic problem of the calculus of variations, namely to find the extremal loop of the following action functional $$F(x):=\int_0^\tau L(t,x,\dot{x})dt,\,\,\,\,\forall x\in W_\tau=W^{1,2}(\mathbf{R}/\tau\mathbf{Z},\mathbf{R}^n),$$ where $L\in C^2((\mathbf{R}/\tau\mathbf{Z})\times \mathbf{R}^n\times \mathbf{R}^n, \mathbf{R})$ such that $L_{p,p}(t,x,p)$ is symmetric and positive definite (Legendre convexity condition). An extreme $x$ of the functional $F$ corresponds to a $\tau$-periodic solution of the Euler-Lagrangian system $$\frac{d}{dt}L_{\dot{x}}(t,x,\dot{x})-L_x(t,x,\dot{x})=0.$$ We define for an extreme $x$, $$P(t):=L_{p,p}(t,x(t),\dot{x}(t)),\,\,\,\,Q(t):=L_{x,p}(t,x(t),\dot{x}(t)),\,\,\,\,R(t):=L_{x,x}(t,x(t),\dot{x}(t)).$$ Then the Hessian of $F$ at $x$ corresponds to the following linear periodic Sturm system $$-(P\dot{y}+Qy)^{\dot{}}+Q^T\dot{y}+Ry=0,$$ which can be switched to a linear Hamiltonian system $\dot{y}=JB(t)y$ with \bea B(t)\equiv B_x(t)=\left(\begin{matrix} P^{-1}(t) & -P^{-1}(t)Q(t) \cr -Q(t)P^{-1}(t) & Q(t)^TP^{-1}(t)Q(t)-R(t)\cr\end{matrix}\right).\nonumber\eea Denote by $\gamma_x$ its fundamental solution. We denote by $m^-(x)$ and $m^0(x)$ respectively the Morse
index and nullity of the functional $F$ at the extreme $x$ in $W_\tau$, which are always finite.

\begin{theorem}(Viterbo\cite{V90}; An-Long\cite{AL98}; \cite{Long02} Chp. 7) Under above notations, we have $$m^-(x)=i_1(\gamma_x),\,\,\,\,m^0(x)=\nu_1(\gamma_x).$$
\end{theorem}

This was later extended to more general relation between Bott index function and Maslov-type index function of $\omega$ for corresponding $\omega$-boundary value problem by Long (see e.g. \cite{Long02}, Chapter 7). The proofs are based on the relation between Duistermaat version of Maslov index and Maslov-type index and their homotopy invariance which can reduce the check to concrete computations for the standard simple cases.

\begin{remark}
One can use also spectral flow to understand the relative Morse index and its relation to Maslov index, please refer to \cite{HS09} and references therein.
\end{remark}

\section{Gutzwiller's Semiclassical Trace Formula}

\subsection{Physical Derivations-WKB Method}

\cite{CAMTV09}, see also \cite{L01}

Feynman's path integrals (\cite{FH65}) provides the most transparent physical derivation of the semiclassical trace formula. The review by Muratore-Ginanneschi (\cite{MG03}) gives a very nice interpretation from this viewpoint and the mathematical theory of functional determinants are used to draw insights in Gutzwiller's trace formula. Here we present another viewpoint (WKB method) more in the spirit of symplectic geometry and Hamiltonian dynamical systems. Both path integral and WKB methods are complementary and useful to get a better understanding on this fundamental formula, and the relation between the two viewpoints is also intriguing.

Recall that the stationary phase principle (SPP) says that the evaluation of the integrals of type $$I=\int d^nx A(x)e^{\sqrt{-1}S(x)/\hbar}$$ is approximated by the contributions from the saddle points of $S(x)$: $$I\approx \sum_{x_k} A(x_k)e^{\sqrt{-1}S(x_k)/\hbar}\frac{(2\pi\sqrt{-1}\hbar)^{n/2}}{\sqrt{\det D^2S(x_k)}}. $$

A quantum mechanical system is characterized by the wave function $\psi(q,t)$ ($q\in\mathbf{R}^n$) determined by the non-relativistic Schr\"{o}dinger equation
$$\sqrt{-1}\hbar\partial_t \psi(q,t)=\hat{H}(q,-\sqrt{-1}\hbar\nabla)\psi(q,t).$$ As is well known that the solutions $\phi_k(q)$($k\in\mathbf{N}$) to the stationary Schr\"{o}dinger equation
$$\hat{H}(q,-\sqrt{-1}\hbar\nabla)\phi_k(q)=E_k\phi_k(q)$$  are orthonormal basis of the complete Hilbert space $L^2(\mathbf{R}^n)$ of states. The wave function $\psi(q,t)$ can be written in terms of these eigenstates as follows (Fourier expansion)
\bes\psi(q,t)&=&\sum_k \left(\int d^nq' \phi_k^*(q')\psi(q',0)\right)e^{-\sqrt{-1}\frac{E_kt}{\hbar}}\phi_k(q)\\
             &=& \int d^nq'\sum_k \phi_k^*(q')\phi_k(q)e^{-\sqrt{-1}\frac{E_kt}{\hbar}}\psi(q',0)\\
             &=& \int d^nq' K(q,q',t)\psi(q',0)
\ees
where $K(q,q',t)=\sum_k \phi_k^*(q')\phi_k(q)e^{-\sqrt{-1}\frac{E_kt}{\hbar}}$ is called the propagator, a fundamental object in quantum physics and is the fundamental solution to the Schr\"{o}dinger equation, i.e., $K(q,q',t)$ solves
$$\sqrt{-1}\hbar\partial_t K(q,q',t)=\hat{H}(q,-\sqrt{-1}\hbar\nabla)K(q,q',t),\,\,\,\,K(q,q',0)=\delta(q-q').$$ So the propagator can be treated as a wave function propagating at $t=0$ from a $\delta$-function at $q=q'$

To relate the propagator to the energy density, one takes the Laplace transform of the propagator $K(q,q',t)$ to get the Green function $G(q,q',E)$, another fundamental object in quantum physics. To evaluate the Green function, one inserts a term with positive $\epsilon$
\bes
G(q,q',E+\sqrt{-1}\epsilon) &=&\frac{1}{\sqrt{-1}\hbar}\int_0^\infty dt e^{\frac{\sqrt{-1}}{\hbar}(E+\sqrt{-1}\epsilon)t}K(q,q',t)\\
                            &=&\frac{1}{\sqrt{-1}\hbar}\int_0^\infty dt e^{\frac{\sqrt{-1}}{\hbar}(E+\sqrt{-1}\epsilon)t}\sum_k \phi_k^*(q')\phi_k(q)e^{-\sqrt{-1}\frac{E_kt}{\hbar}}\\
                            &=& \frac{1}{\sqrt{-1}\hbar}\sum_k \phi_k^*(q')\phi_k(q)\int_0^\infty dt e^{\frac{\sqrt{-1}}{\hbar}(E-E_k+\sqrt{-1}\epsilon)t}\\
                            &=& \sum_k \frac{\phi_k^*(q')\phi_k(q)}{E-E_k+\sqrt{-1}\epsilon}.
\ees
Its trace is
\bes \textrm{Tr} G(q,q',E+\sqrt{-1}\epsilon) &=& \int d^nq G(q,q,E+\sqrt{-1}\epsilon)=\int d^nq \sum_k \frac{\phi_k^*(q')\phi_k(q)}{E-E_k+\sqrt{-1}\epsilon}\\
                                            &=& \sum_k \frac{\delta_{kk}}{E-E_k+\sqrt{-1}\epsilon}=\sum_k \frac{1}{E-E_k+\sqrt{-1}\epsilon}.
\ees
So one can express the density of states in terms of the trace of the Green function for real energies
$$\rho(E)=\sum_k \delta(E-E_k)=-\lim_{\epsilon\rightarrow 0}\frac{1}{\pi}\textrm{Im}\textrm{Tr}G(q,q',E+\sqrt{-1}\epsilon),$$ where we have used the identity $\delta(x-x')=-\lim_{\epsilon\rightarrow 0}\frac{1}{\pi}\textrm{Im}\frac{1}{x-x'+\sqrt{-1}\epsilon}$ ($\epsilon>0$). Hence we are left to evaluate the trace of the Green function to find the energy spectrum.

The strategy here is as follows. One starts from the Wentzel-Kramers-Brillouin(WKB) Ansatz to get the semiclassical wave function via Madlung flow, then gets the semiclassical propagator, namely the semiclassical Van Vlech propagator which is then Laplace transformed into the Green function followed by the first application of the stationary phase principle. We get the final formula by taking the trace of the Green function and applying again the stationary phase principle.

In the semiclassical approximation, the de Broglie wavelength $\lambda\sim \hbar/p$ is short compared with the length scale where the potential varies significantly. This short wave approximation can be developed by taking the limit where $\hbar\rightarrow 0$. Now the wave function can be written in the form of WKB Ansatz: $$\psi_{sc}(q,t)=A(q,t)e^{\sqrt{-1}R(q,t)/\hbar}.$$ Plugging this into the Schr\"{o}dinger equation with natural mechanical Hamiltonian $H(q,p)=\frac{p^2}{2m}-V(q)$ and setting $\rho(t)=A^2(q,t)$ and $v(q,t)=\frac{1}{m}\nabla R(q,t)$, we have in the limit when $\hbar\rightarrow 0$
\bes
\partial_tR &+& H(q,\nabla R(q,t))=0;\\
\partial_t\rho &+& \nabla(\rho v)=0.
\ees
This is the equation of the Madlung flow the first of which is exactly the Hamilton-Jacobi equation! Since $$\frac{d}{dt}R(q(t),t)=\partial_t R(q,t)+\nabla R(q,t)\dot{q}(t)=-H(q,p)+p\dot{q}=L(q(t),\dot{q}(t),t)$$ with $L$ the Lagrangian of the classical system, the evolution of $R(q,t)$ is given by the following integral along the classical trajectory $\gamma$ connecting $q'$ at time $t=0$ to $q$ at time $t$ and solving the corresponding Hamiltonian system $\dot{p}=-\nabla_qH(q,p),\,\,\dot{q}=\nabla_q H(q,p)$
$$ R(q(t),t)=R(q',0)+\int_\gamma dt' L(q(t'),\dot{q}(t'),t):=R(q',0)+R(q,q',t)$$ where $R(q,q',t)$ is Hamilton's principal function. The initial and final momenta are
$$p'=-\nabla_{q'}R(q,q',t),\,\,\,\, p=\nabla_q R(q,q',t).$$ Since the Madlung flow is conserved in time, the infinitesimal volume will be unchanged under the flow, i.e., $$\rho(q(t),t)d^nq=\rho(q',0)d^nq'=\rho(q',0)\det\left(\frac{\partial q'}{\partial q}\right)d^nq.$$ So
\bes \psi_{sc}(q,t) &=& A(q,t)e^{\sqrt{-1}R(q,t)/\hbar}= \sqrt{\det\left(\frac{\partial q'}{\partial q}\right)} A(q',0)e^{\sqrt{-1}(R(q',0)+R(q,q',t))/\hbar}\\
               &=& \sqrt{\det\left(\frac{\partial q'}{\partial q}\right)} e^{\sqrt{-1}R(q,q',t)/\hbar}A(q',0)e^{\sqrt{-1}R(q',0)/\hbar}\\
               &=& \sqrt{\det\left(\frac{\partial q'}{\partial q}\right)} e^{\sqrt{-1}R(q,q',t)/\hbar}\psi_{sc}(q',0)\\
               &=& \sqrt{e^{-\sqrt{-1}\pi i_\gamma(q,q',t)}\left|\det\left(\frac{\partial q'}{\partial q}\right)\right|} e^{\sqrt{-1}R(q,q',t)/\hbar}\psi_{sc}(q',0)\\
               &=& \sqrt{\left|\det\left(\frac{\partial q'}{\partial q}\right)\right|} e^{\sqrt{-1}R(q,q',t)/\hbar-\sqrt{-1}\pi i_\gamma(q,q',t)/2}\psi_{sc}(q',0).
\ees Here $i_\gamma(q,q',t)$ counts the number of the sign changes of the Jacobian determinant along the trajectory $\gamma$ and it is the first hint of some kind of Maslov index! Taking into account the fact that the two points $q$ and $q'$ could be connected by several trajectories $\gamma_k$, the wave function must be summed over all possible paths each of which has its own Maslov index, determinant and Hamilton's principal function, i.e.,
$$\psi_{sc}(q,t)=\sum_{\gamma\,\,\textrm{from}\,\,q'\,\,\textrm{to}\,\,q\,\,\textrm{in\,\,time}\,\,t} \left|\det\left(\frac{\partial q'}{\partial q}\right)\right|^{\frac{1}{2}}e^{\sqrt{-1}R_\gamma(q,q',t)/\hbar-\sqrt{-1}\pi i_\gamma(q,q',t)/2}\psi_{sc}(q',0).$$

We will now derive a semiclassical expression for the propagator $K(q, q', t)$ by considering the
propagator for short time first, and extrapolating from there to arbitrary time $t$.
Recall that $K(q,q',t)$ is the fundamental solution to the Schr\"{o}dinger equation, for infinitesimal short time $\delta t$, away from the singular point $t=0$, we may assume again that \bes K_{sc}(q,q',\delta t)&=& A(q,q',\delta t)e^{\frac{\sqrt{-1}}{\hbar}R(q,q',\delta t)}\\
                                       &\approx& A(q,q',\delta t)e^{\frac{\sqrt{-1}}{\hbar}\left(\frac{m(q-q')^2}{2\delta t}-V(q)\delta t\right)}.
                   \ees
For infinitesimal short time interval, the potential term can be neglected and we are left with a Gaussian distribution with variance $\sigma^2=\frac{\sqrt{-1}\delta t\hbar}{m}$. This can be seen as the finite width approximation to the $\delta$-function since $\delta=\lim_{\sigma\rightarrow 0}\frac{1}{2\pi\sigma^2}e^{-z^2/2\sigma^2}$. This suggests us to take $A=\left(\frac{m}{2\pi\sqrt{-1}\hbar\delta t}\right)^{n/2}$, so we get
$$K_{sc}(q,q',\delta t)\approx \left(\frac{m}{2\pi\sqrt{-1}\hbar\delta t}\right)^{n/2}e^{\frac{\sqrt{-1}}{\hbar}\left(\frac{m(q-q')^2}{2\delta t}-V(q)\delta t\right)}.$$ Because $-\frac{\partial R}{\partial q'}=p'\approx \frac{m(q-q')}{\delta t}, \frac{\partial p'_i}{\partial q_j}=-\frac{\partial^2 R}{\partial q_j\partial q'_i}\approx \frac{m}{\delta t}I$, the factor $\frac{m}{\delta t}$ in $K_{sc}(q,q',\delta t)$ can be interpreted as the determinant of the Jacobian of the transformation from final position coordinates $q$ to initial momentum coordinates $p'$. We now have
\bes K_{sc}(q,q',\delta t)&=& \left(\frac{1}{2\pi\sqrt{-1}\hbar}\right)^{n/2}\left|\det\left(\frac{\partial p'}{\partial q}\right)\right|^{\frac{1}{2}}e^{\sqrt{-1}R(q,q',\delta t)/\hbar}\\
 &=&\left(\frac{1}{2\pi\sqrt{-1}\hbar}\right)^{n/2}\left|\det(-\partial_q\partial_{q'} R_\gamma(q,q',t))\right|^{\frac{1}{2}}e^{\sqrt{-1}R(q,q',\delta t)/\hbar}.
 \ees Using $$K_{sc}(q'',q',t'+\delta t)=\sum_{\gamma\,\,\textrm{from}\,\,q'\,\,\textrm{to}\,\,q''\,\,\textrm{in\,\,time}\,\,t+\delta t} \left|\det\left(\frac{\partial q}{\partial q''}\right)\right|^{1/2}_\gamma e^{\sqrt{-1}R_\gamma(q'',q,t')/\hbar-\sqrt{-1}\pi i_\gamma(q'',q,t')/2}K_{sc}(q,q',\delta t),$$ the additivity of the phases $R(q,q',t)$ and the multiplicativity of Jacobian determinants to evolve our short time approximation of the propagator as in the semiclassical wave function,
we get the final form of the semiclassical propagator, namely the Van Vleck propagator
$$K_{sc}(q,q',t)=\sum_{\gamma\,\,\textrm{from}\,\,q'\,\,\textrm{to}\,\,q\,\,\textrm{in\,\,time}\,\,t}  \frac{1}{(2\pi\sqrt{-1}\hbar)^{n/2}}\left|\det(-\partial_q\partial_{q'} R_\gamma(q,q',t))\right|^{\frac{1}{2}}e^{\sqrt{-1}R_\gamma(q,q',t)/\hbar-\sqrt{-1}\pi i_\gamma(q,q',t)/2}$$ which plays an essential role in the following semiclassical quantization.

Now we derive the Green function from the propagator by Laplace transform which can be done term by term
\bes G(q,q',E) &=& \frac{1}{\sqrt{-1}\hbar}\int_0^\infty dt e^{\frac{\sqrt{-1}}{\hbar}Et}K(q,q',t)= \frac{1}{\sqrt{-1}\hbar}\int_0^\infty dt  e^{\frac{\sqrt{-1}}{\hbar}Et}\sum_\gamma K_\gamma(q,q',t)\\
 &=& \sum_\gamma \frac{1}{\sqrt{-1}\hbar}\int_0^\infty dt e^{\frac{\sqrt{-1}}{\hbar}Et}K_\gamma(q,q',t)\equiv \sum_\gamma G_\gamma(q,q',E).
\ees  We focus on each summand $G_\gamma(q,q',E)$. The integral form suggests us to use the stationary phase principle to
\bes G_\gamma(q,q',E) &=& \frac{1}{\sqrt{-1}\hbar} \int_0^\infty dt e^{\frac{\sqrt{-1}}{\hbar} Et} K_\gamma(q,q',t)\\
&=&  \frac{1}{\sqrt{-1}\hbar} \int_0^\infty dt \frac{\left|\det(-\partial_q\partial_{q'} R_\gamma(q,q',t))\right|^{\frac{1}{2}}}{(2\pi\sqrt{-1}\hbar)^{n/2}}e^{\sqrt{-1}(R_\gamma(q,q',t)+Et)/\hbar-\sqrt{-1}\pi i_\gamma(q,q',t)/2}.
\ees
The stationary condition reads $\partial_t R_\gamma(q,q',t)+E=0$ and the time of the saddle point can be expressed as $t^*=t^*(q,q',E)$ which is just the time that a particle of energy $E$ takes from $q'$ to $q$ by comparing with the Hamilton-Jacobi equation. This implies that the SPP only holds for long time trajectories.  Also we know the pre-exponential contribution $\partial_t^2(R_\gamma(q,q',t)+Et)=\partial_t^2R_\gamma(q,q',t)$.  So, after applying the SPP, we have
\bes G_\gamma(q,q',E)
&=& \sqrt{\frac{2\pi\sqrt{-1}\hbar}{\partial_t^2 R(q,q',t^*)}} \frac{1}{\sqrt{-1}\hbar}  \frac{\left|\det(-\partial_q\partial_{q'} R_\gamma(q,q',t^*))\right|^{\frac{1}{2}}}{(2\pi\sqrt{-1}\hbar)^{n/2}}e^{\sqrt{-1}(R_\gamma(q,q',t^*)+Et^*)/\hbar-\sqrt{-1}\pi i_\gamma(q,q',t^*)/2}\\
&=& \frac{\left|\det(-\partial_q\partial_{q'} R_\gamma(q,q',t^*))\right|^{\frac{1}{2}}}{\sqrt{-1}\hbar (2\pi\sqrt{-1}\hbar)^{(n-1)/2} \sqrt{\partial_t^2 R(q,q',t^*)}} e^{\sqrt{-1}\mathcal{A}_\gamma(q,q',E)/\hbar-\sqrt{-1}\pi i_\gamma(q,q',t^*)/2},
\ees
where $\mathcal{A}_\gamma(q,q',E)$ is the action functional evaluated along the trajectory $\gamma$, i.e., the Legendre transform of the Hamilton's principal function:
$$\mathcal{A}_\gamma(q,q',E)=R_\gamma(q,q',t^*)+Et^*=\int_\gamma dt (L(\gamma(t),\dot{\gamma}(t))+E)=\int_\gamma pdq.$$ One can see that $\nabla \mathcal{A}_\gamma(q,q',E)=\nabla R_\gamma(q,q',t^*)+\partial_{t^*}R_\gamma(q,q',t^*)(\nabla t^*)+E(\nabla t^*)=\nabla R_\gamma(q,q',t^*)$ by the stationary condition. To rewrite the determinant in terms of the action functional, we take the local coordinate around both ends of the trajectory $\gamma$ with the  direction of $q$ being that of the flow. Then $p_\|=\partial {q_\|}R_\gamma(q,q',t)=|p|$  and $p_\perp=\nabla_{q_\perp}R_\gamma(q,q',t)=0$. Putting all these together, we get that $$\det\left(-\partial_q\partial_{q'}R_\gamma\right)=\frac{\partial_t^2R_\gamma}{|\dot{q}||\dot{q}'|}\det\left(-\partial_q\partial_{q'}\mathcal{A}_\gamma\right)\equiv \frac{\partial_t^2R_\gamma}{|\dot{q}||\dot{q}'|}\det D_{\perp,\gamma}(q,q',E),$$ and the semiclassical Green function along $\gamma$
$$G_\gamma(q,q',E)=\frac{|\det D_{\perp,\gamma}(q,q',E)|^{1/2}e^{\sqrt{-1}\mathcal{A}_\gamma(q,q',E)/\hbar-\sqrt{-1}\pi i_\gamma(q,q',E)/2}}{\sqrt{-1}\hbar(2\pi\hbar\sqrt{-1})^{(n-1)/2}\sqrt{|\dot{q}||\dot{q}'|}}.$$

For the trajectory $\gamma$ from $q'$ to $q$ in time $t$ (long trajectory in physics terminology), its contribution to the trace is

\bea \textrm{Tr}G_\gamma(q,q',E)&=&\int d^nq G_\gamma(q,q,E)\nonumber\\
&=& \int d^nq\frac{|\det D_{\perp,\gamma}(q,q,E)|^{1/2}e^{\sqrt{-1}\mathcal{A}_\gamma(q,q,E)/\hbar-\sqrt{-1}\pi i_\gamma(q,q,E)/2}}{\sqrt{-1}\hbar(2\pi\hbar\sqrt{-1})^{(n-1)/2}|\dot{q}|}\nonumber
\eea
The integral form suggests us to use again the stationary phase method. The stationary condition for the exponent $$0=\nabla \mathcal{A}(q,q,E)=\nabla_q \mathcal{A}(q,q',E)|_{q'=q}+\nabla_{q'}\mathcal{A}(q,q',E)|_{q'=q}=p-p'$$ tells us that in fact the starting point and the end point are the same not only in configuration space but also in momentum space! So we get the periodic solutions of the corresponding classical Hamiltonian system. This suggests performing the stationary phase method in the complementary direction of the closed trajectory. It is appropriate to use the local coordinate introduced above. So we get
\bes\textrm{Tr}G_\gamma(q,q',E)&=& \int \frac{dq_\| d^{n-1}q_\perp}{\dot{q}_\|}\frac{|\det D_{\perp,\gamma}(q,q,E)|^{1/2}e^{\sqrt{-1}\mathcal{A}_\gamma(q,q,E)/\hbar-\sqrt{-1}\pi i_\gamma(q,q,E)/2}}{\sqrt{-1}\hbar(2\pi\hbar\sqrt{-1})^{(n-1)/2}}\\
&=&\oint\frac{dq_\|}{\dot{q}_\|}\frac{|\det D_{\perp,\gamma}(q,q,E)|^{1/2}e^{\sqrt{-1}\mathcal{A}_\gamma(q,q,E)/\hbar-\sqrt{-1}\pi i_\gamma(q,q,E)/2}}{|\det D'_{\perp,\gamma}(q,q,E)|^{1/2}\sqrt{-1}\hbar}.\ees Here the action $\mathcal{A}_\gamma(q,q,E)$ depends only on $q_\perp$, in fact only on the orbit. This is due to the definition of the action $A_\gamma(q,q,E)=\oint pdq$ and the extremal condition that $q_\perp=0$. The same holds for the index $i_\gamma(q,q,E)$ which we will discuss in the following remarks. The determinant appearing in the denominator is $\det(\partial^2_{qq}\mathcal{A}+\partial^2_{qq'}\mathcal{A}+\partial^2_{q'q}\mathcal{A}+\partial^2_{q'q'}\mathcal{A})$.

We can also write $|\det D_{\perp,\gamma}(q,q,E)|=|\det(-\partial^2_{q'q}\mathcal{A})|=|\det(\partial^2_{q'q}\mathcal{A})|$ in terms of the monodromy matrix/Poincar\'{e} matrix (or stability matrix in physical terminology) $M_\gamma$. Using the relations $p=\nabla_q \mathcal{A}(q,q',E)$ and $p'=-\nabla_{q'}\mathcal{A}(q,q',E)$, we have
\bes\left(\begin{matrix} \delta q_\perp \cr \delta p_\perp\end{matrix}\right)&=&M_\gamma\left(\begin{matrix} \delta q'_\perp \cr \delta p'_\perp\end{matrix}\right)
=\left(\begin{matrix} M_{\gamma, qq} & M_{\gamma,qp} \cr M_{\gamma,pq} & M_{\gamma,pp}\end{matrix}\right)\left(\begin{matrix} \delta q'_\perp \cr \delta p'_\perp\end{matrix}\right)\\
&=&\left(\begin{matrix} -(\partial^2_{q'q}\mathcal{A})^{-1}(\partial^2_{q'q'}\mathcal{A}) & -(\partial^2_{q'q}\mathcal{A})^{-1} \cr (\partial^2_{qq'}\mathcal{A})-(\partial^2_{qq}\mathcal{A})(\partial^2_{q'q}\mathcal{A})^{-1}(\partial^2_{q'q'}\mathcal{A}) & -(\partial^2_{qq}\mathcal{A})(\partial^2_{q'q}\mathcal{A})^{-1}\end{matrix}\right)\left(\begin{matrix} \delta q'_\perp \cr \delta p'_\perp\end{matrix}\right).
\ees
So
$$\det({I-M_\gamma})=\frac{\det(\partial^2_{qq}\mathcal{A}+\partial^2_{qq'}\mathcal{A}+\partial^2_{q'q}\mathcal{A}+\partial^2_{q'q'}\mathcal{A})}{\det(\partial^2_{q'q}
\mathcal{A})}.$$
Combining these facts together, we get
$$\textrm{Tr}G_\gamma(q,q',E)=\frac{e^{\sqrt{-1}\mathcal{A}_\gamma(E)/\hbar-\sqrt{-1}\pi i_\gamma(E)/2}}{\sqrt{-1}\hbar|\det(I-M_\gamma)|^{1/2}}\oint \frac{dq_{\|}}{\dot{q}_\|}
                                      =\frac{e^{\sqrt{-1}\mathcal{A}_\gamma(E)/\hbar-\sqrt{-1}\pi i_\gamma(E)/2}}{\sqrt{-1}\hbar|\det(I-M_\gamma)|^{1/2}}T_\gamma,$$ with $T_\gamma$ the period of the closed orbit.

The sum should be over all prime closed orbits and their iterations, whence we have totally
$$\textrm{Tr} G(q,q',E)=\textrm{Tr}G_0+\frac{1}{\sqrt{-1}\hbar}\sum_{\gamma: \textrm{closed\,\,orbit}}T_\gamma\frac{e^{\sqrt{-1}\mathcal{A}_\gamma(E)/\hbar-\sqrt{-1}\pi i_\gamma(E)/2}}{|\det(I-M_\gamma)|^{1/2}}$$
with the zero-length contribution ($G_0$-part) to be interpreted as follows.

The SPP cannot be applied when $t^*$ is small because of the divergence of the propagator, so we have to treat this case separately. We evaluate the integral involving the short time form of the semiclassical propagator
\bes G_0(q,q',E)&=&\frac{1}{\sqrt{-1}\hbar}\int_0^\infty dt \left(\frac{m}{2\pi\sqrt{-1}t}\right)^{n/2}e^{\frac{\sqrt{-1}}{\hbar}\left(\frac{m(q-q')^2}{2t}-V(q)t+Et\right)}\\
&=& -\frac{\sqrt{-1}m}{2\hbar^2}\left(\frac{\sqrt{2m(E-V(q))}}{2\pi\hbar|q-q'|}\right)^{\frac{n-2}{2}}H^{(1)}_{\frac{n-2}{2}}(\mathcal{A}_0(q,q',E)/\hbar)
\ees
where $$H_\nu^{(1)}(x):=\frac{1}{\sqrt{-1}\pi}\int_0^\infty dt \frac{e^{\frac{x}{2}(t-\frac{1}{t})}}{t^{\nu+1}}$$ is the Hankel function of the first kind and $\mathcal{A}_0(q,q',E)=\sqrt{2m(E-V(q))}|q-q'|$ is the short distance approximation of the action functional $\mathcal{A}\approx p\Delta q \approx \sqrt{2m(E-V(q))}|q-q'|$. Recall that the Hankel function can be written as $H_\nu^{(1)}(x)=J_\nu(x)+\sqrt{-1}N_\nu(x)$ with $J_\nu(x)$ and $N_\nu(x)$ the Bessel functions of the first and the second kind. The asymptotic form of the Bessel function of the first kind is $J_\nu(x)\approx \frac{1}{\Gamma(\nu+1)}\left(\frac{x}{2}\right)^{\nu}$ for $|x|\ll 1$, and the Bessel function of the second kind $N_\nu(x)$ is singular at the origin which is fortunately not needed here.
\bes
-\frac{1}{\pi}\textrm{ImTr}G_0 &=& -\frac{1}{\pi}\int d^nq \textrm{Im}\lim_{q'\rightarrow q} G_0(q,q',E)\\
&=&\int d^nq \textrm{Im}\lim_{q'\rightarrow q} \frac{\sqrt{-1}m}{2\pi \hbar^2}\left(\frac{\sqrt{2m(E-V(q))}}{2\pi\hbar|q-q'|}\right)^{\frac{n-2}{2}}H^{(1)}_{\frac{n-2}{2}}(\mathcal{A}_0(q,q',E)/\hbar)\\
&=& \int d^nq \lim_{q'\rightarrow q} \frac{m}{2\pi \hbar^2}\left(\frac{\sqrt{2m(E-V(q))}}{2\pi\hbar|q-q'|}\right)^{\frac{n-2}{2}} \textrm{Re} (H^{(1)}_{\frac{n-2}{2}}(\mathcal{A}_0(q,q',E)/\hbar))\\
&=& \int d^nq \lim_{q'\rightarrow q} \frac{m}{2\pi \hbar^2}\left(\frac{\sqrt{2m(E-V(q))}}{2\pi\hbar|q-q'|}\right)^{\frac{n-2}{2}}  J_{\frac{n-2}{2}}(\mathcal{A}_0(q,q',E)/\hbar)\\
&\approx& \int d^nq \lim_{q'\rightarrow q} \frac{m}{2\pi \hbar^2\Gamma(n/2)}\left(\frac{\sqrt{2m(E-V(q))}}{2\pi\hbar|q-q'|}\right)^{\frac{n-2}{2}}\left(\frac{\sqrt{2m(E-V(q)}|q-q'|}{2\hbar}\right)^{\frac{n-2}{2}}\\
&=& \frac{m}{\hbar^n 2^{n-1}\pi^{n/2}\Gamma(n/2)}\int_{V(q)<E} d^nq (2m(E-V(q)))^{\frac{n-2}{2}}.
 \ees
We claim that this is the same as the average density of quantum states for the system. In fact, the density of quantum states is the derivative with respect to the energy of the number $N(E)$ of states under the energy $E$ which can be approximated by the average number $\bar{N}(E)$ which is just the volume of the subset of phase space consisting of states with energy not exceeding $E$, divided by the size of each quantum cell $h^n$ (Weyl rule),
\bes
\bar{N}(E)&=&\frac{1}{h^n}\int d^nqd^np\Theta(E-H(p,q))=\frac{1}{h^n}\int d^nqd^np\Theta(E-\frac{p^2}{2m}-V(q))\\
&=& \frac{1}{h^n}\int_{V(q)<E} d^nq \frac{\pi^{n/2}\sqrt{2m(E-V(q))}^n}{\Gamma(1+n/2)}\\
&=& \frac{\pi^{n/2}}{h^n\frac{n}{2}\Gamma(n/2)}\int_{V(q)<E} d^nq (2m(E-V(q)))^{n/2}
\ees
where $\Theta$ is the Heaviside function and we have used the formula for the volume of sphere $S^{n-1}$ of radius $r$ and dimension $n$ $\textrm{Vol}_{S^{n-1}}=\frac{\pi^{n/2}r^n}{\Gamma(1+n/2)}$ to integrate out the variable $p$ with radius $|p|=\sqrt{2m(E-V(q))}$. Now the average density is
\bes
\frac{d\bar{N}(E)}{dE}&=& \frac{\pi^{n/2}}{h^n\frac{n}{2}\Gamma(n/2)}\int_{V(q)<E} d^nq \frac{n}{2}\cdot 2m\cdot (2m(E-V(q)))^{(n-2)/2}\\
&=& \frac{2m\pi^{n/2}}{h^n\Gamma(n/2)}\int_{V(q)<E} d^nq (2m(E-V(q)))^{(n-2)/2}\\
&=& \frac{m}{\hbar^n 2^{n-1}\pi^{n/2}\Gamma(n/2)}\int_{V(q)<E} d^nq (2m(E-V(q)))^{(n-2)/2},
\ees
as claimed above.

\begin{remarks}
Of course, the biggest problem is how to make this derivation mathematically rigorous, even partially. Maybe there is no hope to achieve this goal, however any attempts in this direction are deserved some of which will be reviewed in the next subsection.

Justification on the energy-time reduction with respect to the index can be found in \cite{HS09}(Lemma 5.3) and references therein. For a natural mechanical system with continuous symmetry (Hamiltonian action of a Lie group), it is an interesting problem to study the interplay between the reduction of the symmetry and the index theory. Discrete symmetries deserve to pay more attention also (see e.g., \cite{FR97}, \cite{R89}).
\end{remarks}

\subsection{Mathematical Justifications}

Semiclassical analysis is the main mathematical tool to justify the semiclassical trace formulas. There are many excellent texts on this topic, please refer to \cite{UW12} for references (see also\cite{G77}).

Colin de Verdi\`{e}re in his thesis (\cite{CdV73a}, \cite{CdV73b}) used the short-time expansion of the Schr\"{o}dinger kernel and a finite dimensional approximation of the Feynman path integral to prove that the spectrum of the Laplacian determines generically the lengths of the closed geodesics. Chazarian (\cite{C74}) derived qualitatively the form of the trace for the wave kernel via Fourier integral operators.

As did in \cite{DG75}, the mathematical derivation via microlocal analysis of the semiclassical trace formula is the applications of the full power of the principal symbol calculus of Fourier integral operators and oscillating integral operators. For a unitary operator $U(t)$, say of the form $e^{-\sqrt{-1}tQ}$ for some first order positive elliptic pseudo-differential operator $Q$ on compact manifold, the trace $\textrm{Tr}(U(t))$ is an oscillating integral composed of the Schwartz kernel $\delta(x,y)$ of the identity and the Schwartz kernel $U(t,x,y)$ of the unitary group with corresponding principal symbols the conormal bundle of the diagonal $\Gamma_\delta=\{(x,\xi;x,\xi)\}\subset T^*(\mathbf{R}^n\times \mathbf{R}^n)$ and the canonical relation belonging to the Hamiltonian flow $F^t$ of $H(x,\xi)$: $\Gamma_U=\{(t,\tau;x,\xi;y,-\eta)\,|\, F^t(y,\eta)=(x,\xi), \tau=-H(x,\xi)\}$. The latter $\Gamma_U$ can be seen as sept out from $\{(0,\tau)\}\times\Gamma_\delta$ by the extended Hamiltonian flow of $\mathcal{H}(t,\tau;x,\xi;y,\eta)=\tau+H(x,\xi)$. Now the principal symbol of the trace under the clean intersection condition is $$\mathcal{P}:=\Gamma_U\circ \Gamma_\delta=\{(t,\tau)\,|\, \exists (x,\xi)\in T^*(\mathbf{R}^n): F^t(x,\xi)=(x,\xi), H(x,\xi)=-\tau\}.$$ Applying the Fourier transform to $\textrm{Tr}(U(t))$, one gets the energy spectral density. This again can be understood as an oscillating integral with principal symbol the canonical transformation $(t, -E)\mapsto (E,t)$ of the principal symbol $\mathcal{P}$ of $\textrm{Tr}(U(t))$. Later this idea was extended and improved by Guillemin-Uribe(\cite{GU89}), Brummelhuis-Uribe (\cite{BU91}), Meinrenken(\cite{M92}), Paul-Uribe(\cite{PU95}) and Combescure-Ralston-Robert (\cite{CRR99}, compare \cite{MW01}) among others via various versions of microlocal analysis.

For recent progress, please refer to \cite{FT15}.

\subsection{Selberg Trace Formula}

\cite{S56}, \cite{H59}, \cite{H87}, \cite{M04}, \cite{H7683}, \cite{BV86}, \cite{CV90}, \cite{Mc72}

The $2$-dimensional surfaces of constant (negative) curvature play an important role in the development of non-Euclidean geometry. It is well known that the free motion (=geodesics) is strongly chaotic which makes it an ideal playground for studying the quantum chaos. On the other hand we know from Riemann surface theory that they are equipped with many discrete symmetries. In 1956, Selberg (\cite{S56}) discovered his formula in the attempt to find a relationship between Riemann's zeta-function and the geometry. The zeros of the zeta-function would play the role of the eigenvalues and the logarithm of the primes are the corresponding periodic orbits. It is the main way to study the fine structure of the spectrum of the Laplacian. It turns out to be a special case of the Gutzwiller trace formula (\cite{H59}). There is a huge mathematical literature devoted to this subject, especially in number theory and harmonic analysis, here we focus on the quantum  manifestations of the chaotic geodesic flows. Gutzwiller (\cite{G80}) has drawn attention to the relation between Selberg's trace formula and the semiclassical expansion of Green's function described by a path integral. He also studied (\cite{G83}) the scattering on a compact surface of constant negative curvature following Lax and Phillips (\cite{LP76}).

Let $\mathbb{H}^2:=\{z=x+\sqrt{-1}y\in \mathbf{C}\,|\, \textrm{Im}(z)>0\}$ be the Poincar\'{e} upper half plane which is a simply connected two-dimensional Riemannian manifold with metric $ds^2=\frac{1}{y^2}(dx^2+dy^2)$ of (constant) Gauss curvature $-1$. The distance $d(z,z')$ between two points $z$ and $z'$ is such that $$\cosh d(z,z')=1+\frac{|z-z'|^2}{2\textrm{Im}(z)\textrm{Im}(z')}.$$ The orientation preserving isometry group of $\mathbb{H}^2$ is $$\textrm{Iso}^+(\mathbb{H}^2)=\textrm{PSL}(2,\mathbf{R})=\textrm{SL}(2,\mathbf{R})/\{\pm I_2\}$$ with action of $\gamma=\left(\begin{matrix} a & b\\ c &d \end{matrix}\right)\in \textrm{Iso}^+(\mathbb{H}^2)$ by the linear fractional transformation (M\"{o}bius transformation) $z\mapsto \gamma\cdot z:=\frac{az+b}{cz+d}$. We call $\gamma\in \textrm{Iso}^+(\mathbb{H}^2)$ hyperbolic if its length defined by $\tau_\gamma=\tau(\gamma):=\inf_{z\in \mathbb{H}^2}d(\gamma z,z)$ is positive, and we call a Fuchsian subgroup $\Gamma\subset \textrm{Iso}^+(\mathbb{H}^2)$ strictly hyperbolic if for any $\gamma\in\Gamma\backslash\{1\}$, $\tau_\gamma>0$. Then the quotient $\Gamma\backslash \mathbb{H}^2$ by a strictly hyperbolic Fuchsian subgroup is a smooth Riemann surface of constant negative curvature, and we denote its fundamental domain by $F$. For hyperbolic element $\gamma\in\Gamma$, it is a fact that the centralizer $\mathcal{Z}_\gamma:=\{g\in \Gamma\,|\, g\gamma=\gamma g\}$ is infinite cyclic subgroup $\mathcal{Z}_\gamma=\{\gamma_*^n\,|\, n\in\mathbf{Z}\}$ where $\gamma_*\in\Gamma$ is the unique element such that $\gamma_*^m=\gamma$ for some $m\in\mathbf{N}$ and no $\tilde\gamma\in\Gamma$ such that $\tilde\gamma^n=\gamma_*$ for any $n\in\mathbf{N}$ and $n>1$. Such $\gamma_*$ is called primitive of $\gamma$. We denote by $H$ the set of conjugacy classes of $\gamma\in\Gamma$ and $H_*\subset H$ the subset of primitive elements.

\begin{theorem}(The Selberg Trace Formula(STF))  Let $h(r)$ be a real function which is analytic on $|\textrm{Im}(r)|\leq \frac{1}{2}+\delta$ such that
$$h(-r)=h(r),\,\,\,\, |h(r)|\leq A(1+|r|)^{-2-\delta},\,\,\,\,(A>0,\delta>0).$$ Then:
$$\sum_{k=0}^\infty h(r_k)=\frac{\textrm{Area}(\Gamma\backslash\mathbb{H}^2)}{4\pi}\int_{-\infty}^\infty h(r)\tanh(\pi r)rdr+\sum_{\gamma\in H_*}\sum_{k=1}^\infty \frac{\tau_\gamma g(k\tau_\gamma)}{2\sinh(k\tau_\gamma/2)},$$ where $r_k=\sqrt{\lambda_k -\frac{1}{4}}$ with $\lambda_k\geq 0$ the spectrum of the Laplacian and $-\frac{\pi}{2}<\arg(r_k)<\frac{\pi}{2}$; $g(t)=\frac{1}{2\pi}\int_{-\infty}^\infty h(r)e^{-\sqrt{-1}rt}dr$ is the Fourier transform of $h$. The sum and the integrals are all absolutely convergent.
\end{theorem}

We start from the elementary derivation of STF in the case where $\Gamma\backslash\mathbb{H}^2$ is a compact Riemann surface following McKean(\cite{Mc72}) and Hejhal(\cite{H87})(see also \cite{BV86}).

The Green function on $\mathbb{H}^2$ is the kernel of the integral operator $(\Delta+\lambda)^{-1}$:
$$((\Delta+\lambda)^{-1}\psi)(z)=\int G(z,z';\lambda)\psi(z')d\nu(z')$$ and can be written as
$$G(z,z';\lambda)=-\frac{1}{2\pi}Q_l(\cosh\tau)=-\frac{1}{2\pi\sqrt{2}}\int_\tau^\infty d\tau' \frac{e^{-\sqrt{-1}\rho\tau'}}{\sqrt{\cosh\tau'-\cosh\tau}}$$
with $Q_l$ the Legendre function of the second kind, $\tau=d(z,z')$, $\lambda=-l(l+1)=\frac{1}{4}+\rho^2$ and $l=-\frac{1}{2}+\sqrt{-1}\rho$. The Green function on $\Gamma\backslash\mathbb{H}^2$ is $$G_F(z,z';\lambda)=\sum_{\gamma\in\Gamma}G(z,\gamma z';\lambda).$$ To assure the convergence of the sum we require that $\textrm{Re} l>0$, or equivalently, $\textrm{Re} \kappa>\frac{1}{2}$ with $\kappa=\sqrt{-1}\rho$ which we use to compensate the exponential proliferation of the orbit under the action of $\Gamma$. The physical meaning of this condition is that we must be deeply in the classically forbidden energy region. When taking the trace of the corresponding integral operator of the Green function, it is still divergent due to the logarithmic singularity at $z'=z$. We regularize this by differentiating with respect to $\lambda$: $$K_0(z,z';\lambda)=-\frac{\partial}{\partial \lambda}G(z,z';\lambda)=-\frac{1}{2\pi}\frac{1}{2\kappa}\frac{\partial}{\partial \kappa} Q_{-1/2+\kappa}(\cosh d(z,z'))\equiv K(\cosh d(z,z')).$$ Now we get the regularized operator $-\frac{\partial}{\partial \lambda}(\Delta+\lambda)^{-1}=(\Delta+\lambda)^{-2}$ with integral kernel $$K_F(z,z';\lambda)=\sum_{\gamma\in\Gamma} K(\cosh d(z,\gamma z')),$$ which we will use to get the trace formula by two ways of evaluation. On the one hand, by using the eigenfunction basis, we have $$\textrm{Tr} K_F=\sum_{k=0}^\infty (\lambda-\lambda_k)^{-2}.$$ On the other hand, $$\textrm{Tr} K_F=\int_F K_F(z,z)d\nu(z)=\sum_{\gamma\in\Gamma} \int_F K(\cosh d(z,\gamma z))d\nu(z).$$ We rewrite it in two steps, the first of which is reorganizing the summation via conjugacy classes.
\bes
\textrm{Tr} K_F&=&\sum_{\gamma\in\Gamma} \int_F K(\cosh d(z,\gamma z))d\nu(z)\\
               &=& \sum_{[\gamma]\in H}\sum_{g\in\Gamma/\mathcal{Z}_\gamma}\int_F K(\cosh d(z,g\gamma g^{-1}z))d\nu(z)\\
               &=& \textrm{Area}(F)\cdot K(1)+\sum_{[\gamma]\in H_*}\sum_{k=1}^\infty \sum_{g\in\Gamma/\mathcal{Z}_\gamma}\int_F K(\cosh d(z,g\gamma^k g^{-1}z))d\nu(z)\\ &=& \textrm{Area}(F)\cdot K(1)+\sum_{[\gamma]\in H_*}\sum_{k=1}^\infty \sum_{g\in\Gamma/\mathcal{Z}_\gamma}\int_{g^{-1}F} K(\cosh d(z,\gamma^kz))d\nu(z)\\
               &=& \textrm{Area}(F)\cdot K(1)+\sum_{[\gamma]\in H_*}\sum_{k=1}^\infty \int_{F_\gamma} K(\cosh d(z,\gamma^kz))d\nu(z),
\ees
where $F_\gamma$ is the fundamental domain for the subgroup $\mathcal{Z}_\gamma$, with respect to which $K(\cosh d(z,\gamma^kz))$ is obviously periodic. Recall also the one to one correspondence between conjugacy classes in the discrete Fuchsian group and closed geodesic on $\Gamma\backslash\mathbb{H}^2$. This step is quite similar to the periodic orbit summations of SCTF (\cite{G71},\cite{BB74}). The second step is evaluating these integrals separately. We choose the half-plane coordinates $\zeta$ to make $\gamma$ diagonal, i.e., reduce it to a dilation $\zeta\mapsto e^{\tau(\gamma)}\zeta$. Note that $\tau(\gamma)$ is the length of the primitive periodic orbit under the action of $\gamma$. Taking the horizontal strip $1<\textrm{Im}\zeta<e^{\tau(\gamma)}$ as the fundamental domain for the group $\mathcal{Z}_\gamma$, we have
\bes
&& \int_{F_\gamma} K(\cosh d(z,\gamma^kz))d\nu(z)\\
 &=& \int_{F_\gamma} \frac{dxdy}{y^2} K(\cosh d(\zeta, e^{k\tau(\gamma)}\zeta))\\
                                               &=& \int_1^{e^{\tau(\gamma)}} \frac{dy}{y^2} \int_{-\infty}^\infty dx K\left(1+\frac{(1-e^{k\tau(\gamma)})^2(x^2+y^2)}{2e^{k\tau(\gamma)}y^2}\right)\\
                                               &=& \tau(\gamma)\int_{-\infty}^\infty dx' K\left(1+2\sinh^2(\frac{k\tau(\gamma)}{2})\cdot (1+x'^2)\right)\,\,\,\,(\textrm{by}\,\, x=yx')\\
                                               &=& \frac{\tau(\gamma)}{\sqrt{2}\sinh(k\tau(\gamma)/2)}\int_{\cosh (k\tau(\gamma))}^\infty \frac{K(t)dt}{\sqrt{t-\cosh (k\tau(\gamma))}}\,\,\,\,(\textrm{by}\,\, t=1+2\sinh^2(\frac{k\tau(\gamma)}{2})\cdot (1+x'^2)).
\ees
By plugging the explicit form of $K(t)$ for the operator $(\Delta+\lambda)^{-2}$:
$$K(\cosh\tau)=-\frac{1}{2\pi\sqrt{2}}\int_\tau^\infty \frac{d\tau'}{\sqrt{\cosh\tau'-\cosh\tau}}\left(\frac{1}{2\kappa}\frac{\partial}{\partial \kappa}\right)e^{-\kappa\tau'},$$ into the above formula, interchanging the order of the integrals and by using $\int_a^b \frac{dt}{\sqrt{(t-a)(b-t)}}=\pi$ with $a=\cosh k\tau(\gamma), b=\cosh\tau', t=\cosh\tau$, we get
$$\int_{F_\gamma} K(\cosh d(z,\gamma^kz))d\nu(z)=\frac{\tau(\gamma)}{2\sinh(k\tau(\gamma)/2)}\left(\frac{1}{2\kappa}\frac{\partial}{\partial \kappa}\right)\frac{e^{-\kappa k\tau(\gamma)}}{2\kappa}.$$ For $K(1)$, we have
\bes
K(1)&=&-\frac{1}{2\pi\sqrt{2}}\int_0^\infty \frac{d\tau'}{\sqrt{\cosh\tau'-1}}\left(\frac{-\tau'}{2\kappa}\right)e^{-\kappa\tau'}=\frac{1}{8\pi\kappa}\int_0^\infty \frac{\tau' e^{-\kappa\tau'}d\tau'}{\sinh(\tau'/2)}\\
&=&\frac{1}{4\pi\kappa}\int_0^\infty\sum_{n=0}^\infty \tau' e^{-(n+\kappa+\frac{1}{2})\tau'} d\tau'=\frac{1}{4\pi\kappa}\sum_{n=0}^\infty (n+\kappa+\frac{1}{2})^{-2}.
\ees
Putting everything together, we get
\bes
\sum_{k=0}^\infty \left(-\frac{\partial}{\partial\lambda}\right)(\lambda-\lambda_k)^{-1}&=&\frac{\textrm{Area}(\Gamma\backslash\mathbb{H}^2)}{2\pi}\left(\frac{1}{2\kappa}\frac{\partial}{\partial \kappa}\right)\frac{\Gamma'(\kappa+\frac{1}{2})}{\Gamma(\kappa+\frac{1}{2})}\\
& &-\sum_{[\gamma]\in H_*}\sum_{k=1}^\infty \frac{\tau(\gamma)}{2\sinh(k\tau(\gamma)/2)}\left(\frac{1}{2\kappa}\frac{\partial}{\partial \kappa}\right)\frac{e^{-\kappa k\tau(\gamma)}}{2\kappa}.
\ees
This is the desired Selberg trace formula for the regularized Green function on the compact Riemann surface $\Gamma\backslash\mathbb{H}^2$ and the quantum mechanical left-hand side is expressed entirely in terms of the area of $\Gamma\backslash\mathbb{H}^2$ and the lengths of its periodic geodesics, hence becomes equal to a classical right-hand side!

However this exact formula is hard to use practically, since the Selberg trace formula is divergent in the vicinity of the eigenvalues. One can get the more general Selberg trace formula by smoothing the trace formula for the Green function via suitable analytic test functions as stated in the theorem. The conditions in the theorem are all needed in the rigorous derivation, and a formal illustrative proof can be found in \cite{BV86} (Appendix L). We note that for Green function, $h(r)=(\lambda-\frac{1}{4}-r^2)^{-1}$ just violates the growth condition in the theorem and that's why we need the regularization procedure above. The most useful choice for $h(r)$ in the general trace formula is the Gaussian $$h(r)=e^{-tr^2},\,\,\,\,g(\tau)=\frac{1}{2\sqrt{\pi}t}e^{-\tau^2/4t}.$$ Now the left-hand side of the Selberg trace formula is the trace of the heat operator $e^{t(\Delta+\frac{1}{4})}$, i.e.,
$$e^{t/4}\sum_{k=0}^\infty e^{-t\lambda_k}=\textrm{Tr} e^{t(\Delta+\frac{1}{4})}=\sum_{k=0}^\infty e^{-t r_k^2},$$ and the trace formula reads
\bes
e^{t/4}\sum_{k=0}^\infty e^{-t\lambda_k}&=& \frac{\textrm{Area}(\Gamma\backslash\mathbb{H}^2)}{(\sqrt{4\pi t})^3}\int_{-\infty}^{\infty}\frac{\tau/2}{\sinh (\tau/2)}e^{-\tau^2/4t}d\tau\\
&&+ \frac{1}{\sqrt{4\pi t}}\sum_{[\gamma]\in H_*}\sum_{k=1}^\infty \frac{\tau(\gamma)}{2\sinh (k\tau(\gamma)/2)}e^{-(k\tau(\gamma))^2/4t}.
\ees
It can be effectively used to study the fine properties of the distribution of eigenvalues, the full asymptotic expansion of the trace for $t\rightarrow 0$ and Weyl's formula by using the classical information. The trace formula can also be used reversely to get the exponential proliferation of periodic orbits by the quantum mechanical ground states of the Laplacian $-\Delta$, i.e., the vanishing of the lowest eigenvalue, an elementary quantum property! For more details, see \cite{BV86} (\S VII and references therein).

For some recent progress on the spectrum of the geodesic flow of negative curvature due to Faure-Tsujii, see \cite{FT15}, \cite{G15} and references therein.

\begin{remark}
The classical Poisson formula can be explained as trace formula. In fact one can use it to prove the trace formula for elliptic curve (i.e. genus $1$ Riemann surface), even for any torus. In this sense Selberg trace formula ia a noncommutative analogue of the Poisson formula (\cite{Mc72}, \cite{M04}) with the role of the lattice in Poisson summation formula replaced by the conjugacy classes of the fundamental group of the Riemann surface.

As observed by Berry-Tabor (\cite{BT76}), the trace formula for the integrable system can be derived from Poisson formula via action-angle coordinates.

\end{remark}

\subsection{Trace Formula and Maslov-type Index}

\cite{M94}

The appearance and anticipation of the Morse index and Maslov index in SCTF is quite earlier.

The WKB ansatz wave function fails at the
turning points of the classical trajectory. A key observation due to Maslov is that while in the Lagrangian formulation of the WKB
ansatz a turning point is singular, switching to the phase space the classical trajectory in
the same neighborhood is smooth. The simplest way to deal with such singularities is as follows: follow the classical trajectory in
$q$-representation until the WKB approximation fails close to the turning point; then insert the transform from coordinate representation to the momentum representation near the turning point and follow the classical trajectory in the $p$-space until encountering
the next $p$-space turning point; go back to the $q$-space representation, an so on.
One can evaluate the transform to the leading order of $\hbar$ by the stationary phase method which involves a Fresnel integral $$\frac{1}{\sqrt{2\pi}}\int_{-\infty}^{+\infty}dx e^{-\frac{x^2}{2\sqrt{-1}a}}=\sqrt{\sqrt{-1}a}=|a|^{1/2}e^{\sqrt{-1}\frac{\pi}{4}\frac{a}{|a|}},$$ yielding an extra $e^{\sqrt{-1}\pi/4}$ phase shift. The single-valuedness of the wave-function leads to the WKB quantization, including the Bohr-Sommerfeld
quantization condition.

Mathematically, they appear gradually in the works of Colin de Verdi\`{e}re, Duistermaat-Guillemin, Duistermaat, Meinrenken and many others.

We recall here the proof due to Meinrenken on the fact that the right version of Maslov index is Conley-Zehnder or Maslov-type index. Here we are content with the restriction on the $\mathbf{R}^{2n}=T^*(\mathbf{R}^n)$ to give the idea, the generalization of which will be commented and addressed elsewhere.

We start with H\"{o}mander's construction of Maslov's principal bundle $\mathcal{M}$ over the Lagrangian Grassmannian $\Lambda(n)$. Fix $L_1\in \Lambda(n)$, by Proposition \ref{signature}, $$\frac{1}{2}(s(L_1,L_2,M_1)-s(L_1,L_2,M_2))=\frac{1}{2}(s(M_1,M_2,L_2)-s(M_1,M_2,L_1))$$ are locally constant and integer-valued when $L_i$ and $M_j$ are transversal which can be used as transition functions to define the principal bundle $\mathcal{M}$. A section of $\mathcal{M}$ over an open subset $U\subset \Lambda(n)$ can be regarded as a function $\phi: U\times \Lambda(n)\rightarrow \frac{1}{2}\mathbf{Z}$ such that $\phi(L_2,M)-\frac{1}{2}s(L_1,L_2,M)$ is independent of $M$ and $\phi(L_2,M)$ is continuous on the set of $M$ which is transversal to $L_i$. Each point $L_2^0\in\Lambda(n)$ determines the germ of a trivialization. One can see this as follows: letting $U$ be some contractible neighborhood of $L_2^0$, and taking $\phi(L_2^0,M)=\frac{1}{2}s(L_1,L_2^0,M)$ and parallel transport on $\mathcal{M}$ yield a canonical trivialization $$\phi(L_2,M)=\frac{1}{2}s(L_1,L_2,M)+[L_1:L_2(t)]$$ with $L_2(t)$ any path in $U$ from $L_2^0$ to $L_2$. Maslov's principal bundle for a Lagrangian submanifold $N$ of $T^*\mathbf{R}^n$ is defined similarly. Let $L_1$ be the vertical polarization and $L_2$ the tangent bundle of $N$. Let $X_i=\mathbf{R}^{n_i}$ ($i=1,2,3$) for simplicity which can be extended to any manifold without difficulty. Suppose $N_2,N_1$ be Lagrangian submanifolds of $T^*(X_3\times X_2)$ and $T^*(X_2\times X_1)$ with sections $\phi_i$ and their Maslov bundles $\mathcal{M}_i$ respectively. Let $S\subset T^*(X_2\times X_2)$ be the conormal bundle of the diagonal $\Delta$. Then $N_2\circ N_1$ is, by definition, the image of $(T^*X_3\times S\times T^*X_1)\cap (N_2\times N_1)$ under the symplectic reduction $\rho: T^*X_3\times S\times T^*X_1 \rightarrow T^*(X_3\times X_1)$ which is an immersed Lagrangian manifold if the intersection is clean. The composed section $\phi_2\circ \phi_1$ is defined to be $$\phi_2\circ \phi_1(W_p)=(\phi_2\times \phi_1)((T_z\rho)^{-1}(W_p))$$ for $W_p\in \Lambda(T_p(T^*(X_3\times X_1)))$ and arbitrary $z\in \rho^{-1}(p)$.

Denote by $V^\mathbf{R}$ and $V^{\mathbf{R}^n}$ the vertical polarization in $T^*\mathbf{R}$ and $T^*(\mathbf{R}^n)$ respectively. We also use $V$ to denote general vertical polarization. The canonical trivialization of the Maslov bundle over $\Gamma_\delta$ is defined to be $$\phi_\delta(W_z)=\frac{1}{2}s(V_z,T_z\Gamma_\delta,W_z).$$ Parallel transport along the solution curve $\tilde\gamma$ of $X_\mathcal{H}$ induces a trivialization of the Maslov bundle over $\Gamma_U$:
$$\phi_U(W_{\tilde\gamma(T)})=[V_{\tilde\gamma(t)}:T_{\tilde\gamma(t)}\Gamma_U]+\frac{1}{2}s(V_{\tilde\gamma(T)},T_{\tilde\gamma(T)}\Gamma_U,W_{\tilde\gamma(T)}).$$
By Proposition \ref{index for pairs} (6), one can see that the first term is just $[V^{\mathbf{R}^n}: TF^t(V^{\mathbf{R}^n})]$. The second term is exactly the canonical trivialization described above when $\Gamma_U$ is transversal to the vertical trivialization.

Suppose that $(T,-E)\in \mathcal{P}$ is a nondegenerate periodic solution $\gamma$ with $T\neq 0$. The Maslov phase $\sigma_\gamma$ appearing in the trace formula is the transition function from the composed section $\phi_U \circ \phi_\delta$ of the Maslov bundle to the canonical trivialization at $(E,T)$. Note that on the symbol level, taking Fourier transform in the final step to get the trace formula corresponds to switch from the horizontal polarization to vertical polarization of $T^*\mathbf{R}$. Evaluating the composed section on $Z:=\textrm{span}(\partial/\partial t)$, we get
\bes
              &&(\phi_U\circ \phi_\delta)(Z) \\
              &=&(\phi_U\times \phi_\delta)(Z\times TS)\\
              &=& [V^{\mathbf{R}^n}:TF^t(V^{\mathbf{R}^n})]+\frac{1}{2}s(V^\mathbf{R}\times V, T\Gamma_U\times T\Delta, Z\times TS)\\
              &=& [V^{\mathbf{R}^n}:TF^t(V^{\mathbf{R}^n})]+\frac{1}{2}s(V^\mathbf{R}\times V, T\Gamma_U\times T\Delta,V^\mathbf{R}\times TS)+\frac{1}{2}s(V^\mathbf{R},Z,T\mathcal{P})\,\,(\textrm{Proposition} \ref{signature})\\
              &=& [V^{\mathbf{R}^n}:TF^t(V^{\mathbf{R}^n})]+\frac{1}{2}s(V^\mathbf{R}\times V, T\Gamma_U\times T\Delta,V^\mathbf{R}\times TS)+\frac{1}{2}\textrm{sgn}(\partial T/\partial E)\\
              &=& [V^{\mathbf{R}^n}:TF^t(V^{\mathbf{R}^n})]+\frac{1}{2}s(V, \Gamma(TF^T), T\Delta)+\frac{1}{2}\textrm{sgn}(\partial T/\partial E)\,\,(\textrm{Proposition} \ref{signature}).
\ees
After applying the final step Fourier transform, $$\sigma_\gamma =[V^{\mathbf{R}^n}:TF^t(V^{\mathbf{R}^n})]+\frac{1}{2}s(V, \Gamma(TF^T), T\Delta).$$
This is exactly the Conley-Zehnder index $i_1(\gamma)$ by taking $V=V^{\mathbf{R}^n}\times V^{\mathbf{R}^n}$ (compare (\ref{CZ=Maslovpahse}))!

\section{Conclusions and Speculations}

The physical and mathematical derivations of Gutzwiller's semiclassical trace formula are notoriously difficult, and many recent advances in mathematics and quantum physics draw inspirations from the attempting to understand this marvelous formula.

We address ourselves here several problems which seem promising.

There are many works devoted to the computations of the Maslov indices for concrete examples like harmonic oscillators, H\'{e}non-Heiles, we will use Maslov-type index developed by Long to consider these examples in more details (monodromy matrices, basic normal forms, summation over the iterations of the periodic orbits) in a forthcoming paper to get a better understanding on the energy spectrum. An intermediate step is to understand the relationships among the various versions of the definitions of Maslov indices as did in the previous \cite{CLM94} with the recent progress in mind.   

It is a common wisdom that Gutzwiller's trace formula does not converge in any sense. Berry (\cite{Be91}) has emphasized this point with the suggestion of the use of the resurgence theory \`{a} la \'{E}calle (see \cite{S14} and references therein). This tool was later used by Berry-Howls (\cite{BH94}) to study the quantum billiards. In an ongoing project, we try to combine the full power of the resurgence theory and the calculations of Maslov-type index theory to get better understanding on the contributions of the periodic orbits to the distribution of the spectrum. Pad\'{e} approximant is also used for this purpose (see e.g. \cite{MaWu01}).

Although the natural setting for the semiclassical trace formula in physics is the cotangent bundle of the configuration space, it seems plausible to develop the Maslov-type index theory especially the iteration theory on general symplectic manifold. Meinrenken has made some progress in this direction. P. Seidel's graded Lagrangian manifold (\cite{S00}, inspired by M. Kontsevich) will be useful, and this is a conceptually intrinsic formulation to play with the fundamental group of the Lagrangian Grassmannian.

Why is Selberg's trace formula exact? Maybe it is due to the higher degree of group symmetry of the Riemann surface, however the precise connection of these two aspects seems inconclusive. It deserves to pursue this further. More general question would be when the semiclassical trace formula is exact.

We use several times of the stationary phase principle in the derivation of the trace formula, and this is one reason that we get finally an asymptotic formula. Recent investigation about the supersymmetry in quantum mechanics and quantum field theory in general leads often to exact stationary approximation. This suggests that in supersymmetric quantum mechanics anticipation of exact trace formula is reasonable, and please refer to \cite{LKK99} for the literature. The idea about the complexification of actions and measures in Feynman path integrals seems promising (\cite{BDSSU15}).

Many natural mechanical systems have symmetries like the rotational invariance of the Hamiltonian under the action of the Lie group $SO(3)$, so the periodic orbits are degenerate and appear in families. In quantum physics, Faddeev-Popov (\cite{FP67}) and Batalin-Vilkovisky (\cite{BV81}) developed the methods to deal with the degenerate Lagrangians with respect to Lie group action in quantizing the gauge theory. Schwarz (\cite{Sch93}) developed the semiclassical analysis in Batalin-Vilkovsky formalism. Very recently, on mathematical side, shifted sympletic structures are developed by Pantev-To\"{e}n-Vaqui\'{e}-Vezzosi (\cite{PTVV13}). This would shed light on the rigorous derivation of the Gutzwiller's semiclassical trace formula.

\end{document}